\documentclass[11pt]{article}
\pdfoutput=1
\usepackage[hmargin=1.4in]{geometry}
\usepackage[utf8]{inputenc}
\usepackage[T1]{fontenc}
\usepackage{Baskervaldx,bbold}
\usepackage[shortlabels]{enumitem}

\usepackage{amsmath,amssymb,amsthm,bm}
\usepackage{mathtools}
\usepackage{quiver}

\newcommand{\zigzag}{\leftrightsquigarrow}

\renewcommand{\to}{\mathchoice{\longrightarrow}{\rightarrow}{\rightarrow}{\rightarrow}}
\newcommand{\from}{\mathchoice{\longleftarrow}{\leftarrow}{\leftarrow}{\leftarrow}}
\let\oldmapsto\mapsto
\renewcommand{\mapsto}{\mathchoice{\longmapsto}{\oldmapsto}{\oldmapsto}{\oldmapsto}}
\newcommand{\To}{\mathchoice{\Longrightarrow}{\Rightarrow}{\Rightarrow}{\Rightarrow}}

\newcommand{\op}{\mathrm{op}}
\newcommand{\pseudo}{\mathrm{ps}}
\newcommand{\strict}{\mathrm{strict}}
\newcommand{\cat}[1]{\mathsf{#1}}
\newcommand{\twocat}[1]{\mathcal{#1}}
\newcommand{\CAT}[1]{\mathsf{#1}}
\newcommand{\Set}{\CAT{Set}}
\newcommand{\Cat}{\CAT{Cat}}
\newcommand{\Diag}{\operatorname{\CAT{Diag}}_{\rightarrow}}
\newcommand{\DiagOp}{\operatorname{\CAT{Diag}}_{\leftarrow}}
\newcommand{\strongDiag}{\operatorname{\CAT{Diag}}_{\rightarrow}^\pseudo}
\newcommand{\strongDiagOp}{\operatorname{\CAT{Diag}}_{\leftarrow}^\pseudo}
\newcommand{\strictDiag}{\operatorname{\CAT{Diag}}_{\rightarrow}^\strict}
\newcommand{\strictDiagOp}{\operatorname{\CAT{Diag}}_{\leftarrow}^\strict}
\newcommand{\Lift}[2]{\operatorname{Lift}^{#2}#1} % lifts of #1 along #2
\newcommand{\iclass}{\cat{N}}
\newcommand{\pclass}{\cat{W}_{\rightarrow}}
\newcommand{\eosclass}{\cat{W}_{\leftarrow}}
\newcommand{\enrichedcat}[1]{#1\text{-}\Cat}
\newcommand{\algebra}[1]{#1\text{-}\mathcal{A}\textit{lg}}

\DeclareMathOperator{\cod}{cod}

\DeclareMathOperator{\id}{id}
\DeclareMathOperator*{\colim}{colim}

\usepackage[hidelinks,bookmarksnumbered]{hyperref}
\usepackage{subcaption}
\usepackage[capitalize,noabbrev,nameinlink]{cleveref}

\newtheorem{theorem}{Theorem}[section]
\newtheorem{proposition}[theorem]{Proposition}
\newtheorem{lemma}[theorem]{Lemma}
\newtheorem{corollary}[theorem]{Corollary}
\theoremstyle{definition}
\newtheorem{definition}[theorem]{Definition}
\newtheorem{example}[theorem]{Example}
\newtheorem{remark}[theorem]{Remark}
\usepackage{etoolbox}
\AtEndEnvironment{definition}{\null\hfill$\lrcorner$}
\AtEndEnvironment{example}{\null\hfill$\lrcorner$}
\AtEndEnvironment{remark}{\null\hfill$\lrcorner$}

\usepackage[style=alphabetic,maxnames=99,maxalphanames=4,maxbibnames=10]{biblatex}
\renewbibmacro{in:}{}

\usepackage{tikz}
\usepackage{tikz-cd}
\usetikzlibrary{arrows}

\title{The diagrammatic presentation of equations in categories}
\author{Kevin Arlin \and James Fairbanks \and Tim Hosgood \and Evan Patterson}
\date{17 January 2024}

\usepackage{tocloft}
\begin{document}

\maketitle

\begin{abstract}
  A lift of a diagram $D\colon\mathsf{J}\to\mathsf{X}$ in a category against a discrete opfibration $\pi\colon\mathsf{E}\to\mathsf{X}$ can be interpreted as presenting a solution to a system of equations in $\mathsf{X}$ presented by $D.$ With this interpretation in mind, it is natural to ask if there is a notion of equivalence of diagrams $D\simeq D'$ that precisely captures the idea of the two diagrams ``having the same solutions''. We give such a definition, and then show how the localisation of the category of all diagrams in $\mathsf{X}$ along such equivalences is isomorphic to the localisation of the slice category $\mathsf{Cat}/\mathsf{X}$ along the class of initial functors. Finally, we extend this result to the 2-categorical setting, proving the analogous statement for any locally presentable 2-category in place of $\mathsf{Cat}$.
\end{abstract}

\tableofcontents

\section{Introduction}
\label{section:introduction}

A diagram in a category can be viewed as a presentation of a system
of equations given by equating parallel paths in its image.
For example, we can present the equation
\[
	y^2 = x^3+x+1
\]
describing an elliptic curve in an affine plane $\mathbb{A}^2$ over a field $k$ by the diagram
% https://q.uiver.app/#q=WzAsMixbMCwwLCJcXG1hdGhiYiBBXjIiXSxbMSwwLCJrIl0sWzAsMSwieV4yIiwwLHsiY3VydmUiOi0yfV0sWzAsMSwieF4zK3grMSIsMix7ImN1cnZlIjoyfV1d
\[\begin{tikzcd}
	{\mathbb A^2} & k
	\arrow["{y^2}", curve={height=-12pt}, from=1-1, to=1-2]
	\arrow["{x^3+x+1}"', curve={height=12pt}, from=1-1, to=1-2]
\end{tikzcd}.\]
What is crucial here is that this diagram \emph{does not commute}.
The locus on which the equation $y^2=x^3+x+1$ holds, which for a geometer
is the variety determined by the system, is given
categorically by the equalizer; in general, the limit of a
diagram in a category may be viewed as the object of all solutions to the
system of equations the diagram presents.

From previous work \cite{patterson2022}, we are motivated by partial differential equations, a field in
which mathematicians must in general study only points, or families of points, in the solution space.
An individual solution is given in familiar categories by a
cone over the diagram tipped by a terminal object or monoidal unit. More generally,
a family of solutions indexed by an arbitrary object $x$ is given
by a cone over the diagram tipped by $x$. Thus we might, for instance, obtain
a two-parameter family of solutions to a system of differential
equations using a cone tipped by $\mathbb R^2$.

The categorical framing allows us to generalize yet further. Letting $D: \cat{J} \to \cat{X}$ be a diagram in a category $\cat{X}$, we shift perspective to consider the choice of a cone over $D$ tipped by $x$
as a \emph{lift} of the functor $D$ into the coslice category $x/\cat{X}$:
% https://q.uiver.app/#q=WzAsMyxbMCwxLCJcXHNxdWFyZSJdLFsxLDEsIlxcY2F0IFgiXSxbMSwwLCJ4L1xcY2F0IFgiXSxbMiwxXSxbMCwxLCJEIiwyXSxbMCwyLCIiLDAseyJzdHlsZSI6eyJib2R5Ijp7Im5hbWUiOiJkb3R0ZWQifX19XV0=
\[\begin{tikzcd}
	& {x/\cat{X}} \\
	\cat{J} & {\cat{X}}
	\arrow[from=1-2, to=2-2]
	\arrow["D"', from=2-1, to=2-2]
	\arrow[dotted, from=2-1, to=1-2]
\end{tikzcd}.\]

We are led to consider generalized solutions via lifting against
functors into $\cat{X}.$ For instance, if the spaces $\mathbb R^n$
are in our category $\cat{X}$ then, by considering lifts against
the canonical functor $\coprod_n (\mathbb R^n/\cat{X})\to \cat{X}$, we can reify
the concept of ``an $n$-parameter family of solutions for some $n$.'' If the interval $[0,1]$
is available in $\cat{X}$ then, by the same token, by lifting against the projection $C\to X$,
where $C$ is the coequalizer of the canonical functors
$\mathbb R^1\times [0,1]/\cat{X}\to \mathbb R^1/\cat{X}$, we reify the concept of
``a homotopy equivalence class of one-parameter families of
solutions.'' It is thus rather natural to close the
family of coslice projections $x/\cat{X}\to \cat{X}$ under arbitrary
colimits and consider notions of solution resulting from
lifting against all of the resulting functors into $\cat{X}$.
The result of this colimit closure is precisely the class of discrete opfibrations
over $\cat{X}$. (We recall the definition of discrete opfibration
below in \cref{definition:discrete-opfibration}.) For reasons that will become clear following \cref{definition:equivalence-of-systems} below, it is infelicitous to continue
generalizing to lifts against arbitrary functors over $\cat{X}$, and so we have found our basic setting.
This was already the setting that emerged from a recent study in applied category
theory by Patterson, Baas, Hosgood, and Fairbanks \cite{patterson2022}.

To summarize thus far: we are interested in lifts
of diagrams $D\colon\cat J\to\cat{X}$ against discrete opfibrations
$\pi\colon\cat E\to \cat{X}$. We take as our central question
under what circumstances such a diagram $D$
admits ``the same'' lifts against all such discrete opfibrations---in equational terms,
the same solutions for every way of interpreting the variables when the diagram is
viewed as a system of equations---as another diagram $D'\colon\cat J'\to \cat{X}.$
Since we want to allow the shape of the diagram to vary,
we are led to consider the category $\DiagOp(\cat{X})$
of all diagrams in $\cat{X}$. (See \cref{section:diagram-categories} for the details on
this category; note for now that it is not simply the slice $\Cat/\cat{X}$ but incorporates
a lax aspect in its morphisms.)
The category $\DiagOp(\cat{X})$, as a suitable domain
for a global limit functor, goes all the way back to
Eilenberg and Mac Lane \cite{eilenberg1945}. It is nevertheless not very
extensively studied (though see \cite{guitart1977} for a characterization of $\DiagOp(\cat{X})$ as a lax cocompletion and \cite{mesiti2023} for
a description of colimits in $\DiagOp(\cat{X})$),
and the notion of limit is more often formalized
``locally,'' one diagram shape at a time. The diagram category has nonetheless arisen repeatedly
in the short history of applied category theory thus far, both
in the current context and in the modeling of functorial
data migrations \cite{spivak2023}; in both cases at the heart of the usefulness
of $\DiagOp(\cat{X})$ is the fact that a diagram in $\cat{X}$ gives a convenient presentation of a
discrete opfibration over $\cat{X}$---for computational purposes, importantly, a presentation that
might well be finite even when $\cat{X}$ is not.

Any morphism $D\to D'$ in $\DiagOp(\cat{X})$ induces a function from
lifts of $D$ against a discrete opfibration $\pi$ to lifts of $D'$ against $\pi$, and we
can use this canonical function to define a class of
\emph{weak equivalences} on $\DiagOp(\cat{X})$ as those morphisms which induce
bijections of lifts against arbitrary discrete opfibrations over $\cat{X}$.
Such maps will be viewed as isomorphisms of the systems of equations presented by the two diagrams.
A very similar class of weak equivalences was studied in Section 9 of
\cite{patterson2022}.

The primary work of this paper, leading to \cref{theorem:localisation-of-slice-cat-and-of-diag-op},
is to compute the localization $\DiagOp(\cat{X})[\cat W^{-1}]$, where $\cat W$ is the class of weak
equivalences described above. The \emph{localization} (\cref{definition:localisation}) of a
category at some class of morphisms universally adjoins inverses to those morphisms \cite{gabriel1967}.
In particular, the localization $\DiagOp(\cat{X})[\cat W^{-1}]$, universally turning the equivalences of systems of equations
into isomorphisms, can be thought of as ``the category of systems of equations in $\cat{X}$.''
From a theoretical point of view, we obtain a completely satisfying answer:

\begin{quote}
\textit{The category $\DiagOp(\cat{X})[\cat W^{-1}]$ of systems of equations in a category
$\cat{X}$ is equivalent to the full subcategory of the slice $\Cat/\cat{X}$ spanned by the discrete opfibrations.}
\end{quote}

It was not at all obvious a priori to these authors that this theorem should hold,
although it is reads quite naturally after the fact.
In particular, the proof includes some technical work in establishing that the morphisms in the localization
can all be represented by \emph{strict} maps in $\Cat/\cat{X}.$

So, the motto is that a ``a system of equations in $\cat{X}$ is a discrete opfibration over $\cat{X}$."
That said, we are especially interested in computationally tractable descriptions of this
localization for future work in algorithmic rewriting for PDEs.
The description of the localization in terms of discrete opfibrations (or, equivalently, copresheaves)
is not appropriate for this purpose, due to the loss of the opportunity for finite
presentation of objects mentioned above.
We thus expect to find more applications for the description of the localization as the category
whose objects are diagrams in $\cat{X}$ and whose morphisms $D\to D''$ are the zigzags
$D\leftarrow D'\to D''$ with the reversed arrow determined by an initial functor
between the domains of $D$ and of $D'.$
Initial functors (\cref{definition:initial-morphism}) can be described either as functors,
restriction along which does not change the limit of a diagram,
or as those left orthogonal to discrete opfibrations.
The other key aspect of the proof of our main theorem,
in addition to the strictification mentioned above,
is to show that inverting diagram morphisms given by initial functors suffices to invert all
weak equivalences---even though not nearly all weak equivalences are themselves given by
initial morphisms.

Our first main result, \cref{theorem:localisation-of-slice-cat-and-of-diag-op}, is not sufficiently general
to apply to most systems of interest in physics and engineering, including those in \cite{patterson2022}.
The reason is that diagrams in a category without extra structure can only describe equations
involving unary operations, whereas equations in applied mathematics tend to involve
operations or differential operators of higher arity. Thus, one often prefers to
take diagrams in a categorical structure going beyond a bare category.

For instance, as simple an equation as $B(x,y)=0$, where $B$ is a bilinear form on a vector space $V$, is most naturally modelled by considering the symmetric monoidal or multicategorical structure on vector spaces. We generally prefer the latter, though we address both, so that we would draw this equation diagrammatically like as
% https://q.uiver.app/#q=WzAsMixbMCwwLCIoVixWKSJdLFsxLDAsImsiXSxbMCwxLCJCIiwyLHsiY3VydmUiOjJ9XSxbMCwxLCIwIiwwLHsiY3VydmUiOi0yfV1d
\[\begin{tikzcd}
	{(V,V)} & k
	\arrow["B"', curve={height=12pt}, from=1-1, to=1-2]
	\arrow["0", curve={height=-12pt}, from=1-1, to=1-2]
\end{tikzcd}.\]
For various types of equations, we might prefer to work in multicategories, perhaps symmetric or
cartesian or equipped with coproducts, or the analogous monoidal,
symmetrical monoidal, cartesian monoidal, or distributive monoidal categories.
Motivated by these many possibilities, the final
\cref{section:2-categorical-story} generalizes our main result from the
2-category of categories to an arbitrary locally presentable 2-category. This
setting directly includes the multicategorical examples and
indirectly includes the others, as we will explain later using standard tricks from 2-monad theory.

To even state the broadened result, the concepts of discrete opfibration and
initial functor must be defined in a general 2-category. These concepts sometimes
generalize in more than one useful way; for our purposes, we propose a definition
of discrete opfibration in an arbitrary 2-category that is inspired by Riehl and
Verity's treatment of formal \mbox{($\infty$-)category} theory in a virtual equipment \cite{riehl-verity}.
This definition reduces to the standard one in the 2-category of categories. We then \emph{define} the initial morphisms to be the class of
morphisms left orthogonal to discrete opfibrations. Thus, what was a theorem for
categories, the existence of Street and Walters' \emph{comprehensive factorization} \cite{street1973},
becomes essentially a definition, using standard features of factorization systems in locally presentable categories.

With these definitions and a limited amount of groundwork, the proofs in
\cref{section:localisations-of-categories-of-diagrams} go through essentially
without change to establish a wide generalization of our main theorem:
\begin{quote}
  \textit{For any locally presentable 2-category $\twocat{K}$ containing an object $x$, the category of discrete opfibrations over $x$, as a full subcategory of the slice $\twocat K/x$, is equivalent to the localization of the diagram category $\DiagOp(x)$ at a class of weak equivalences defined just as in the case $\twocat K=\Cat.$}
\end{quote}

The results of this paper illustrate how relatively deep categorical mathematics, including localizations, 2-dimensional monad theory, and comprehensive factorization systems, can be used to give a novel notion of sameness for systems of equations. In
future work, we intend to make use of these results to implement algorithms for the simplification of such systems
that would be essentially indescribable without the categorical formalism.

\paragraph{Conventions} We write categories in sans-serif, e.g. $\cat{X}$ and $\cat{J}$, functors and objects of an arbitrary 2-category
in serif capitals, e.g. $D$ and $X$, morphisms of a category or a 2-category in serif lower case, e.g. $f$ and $g$, and natural transformations
and other 2-morphisms in Greek script.

While we would philosophically lean toward diagrammatic order of composition, we cannot bring ourselves to apply functions or functors
on the left. So, in the $\Cat$-centric first four sections of the paper, the
juxtaposition $FG$ means to apply $F$ \emph{after} $G$, in the classical order.
In \cref{section:2-categorical-story}, however, we will be working in an abstract 2-category with no
ability to apply 1-morphisms to elements of their domain, and so we will use the diagrammatic order of composition. We highlight
the intent to compose diagrammatically with $\cdot$ but sometimes suppress the $\cdot$ in more complicated composites when
no ambiguity can arise. We make no notational distinction between horizontal and vertical composition
of 2-morphisms, taking the perspective that both are just special cases of the fundamental composition operation---that of pasting---on
a 2-category.

\paragraph{Acknowledgments}

Arlin, Fairbanks, and Patterson acknowledge support from DARPA Award
HR00112220038. Hosgood thanks Evan Cavallo and Ivan di Liberti for helpful
conversations.

\section{Diagram categories}
\label{section:diagram-categories}

Let $\cat{X}$ be an arbitrary category.

\begin{definition}
\label{definition:diag-and-diagop}
  The (oplax) \emph{category of diagrams in $\cat{X}$}, or \emph{diagram category}, denoted $\Diag(\cat{X})$, is defined as follows:
  \begin{itemize}
    \item the objects are diagrams $D\colon\cat{J}\to\cat{X}$ indexed
    by a small category $\cat{J}$, which we denote by $(\cat{J},D)$ (or, more often, just $D$ when $\cat{J}$ is implicit);
    \item the morphisms $(\cat{J},D)\to(\cat{K},E)$ are pairs $(R,\rho)$, where $R\colon\cat{J}\to\cat{K}$ is a functor and $\rho\colon D\To ER$ is a natural transformation.
    (See Figure 1.)
  \end{itemize}
  Dually, we define the \emph{contravariant category of diagrams} $\DiagOp(\cat{X})$ as $\Diag(\cat{X}^\op).$ Thus $\DiagOp(\cat{X})$ has the same objects as
  $\Diag(\cat{X})$, but have morphisms $(\cat{J},D)\to(\cat{K},E)$ given by
  pairs $(R,\rho)$, where $R\colon\cat{K}\to\cat{J}$ is a functor and $\rho\colon DR\To E$ is a natural transformation.
  (See \cref{fig:morphisms-in-diag-and-diag-op}).

  We say that a morphism $(R,\rho)$ in $\Diag(\cat{X})$ or $\DiagOp(\cat{X})$ is \emph{pseudo}\footnote{One could also use the terminology ``\emph{strong}'' but we opt for convention most common in the $2$-categorical literature.} if $\rho$ is a natural isomorphism and is \emph{strict} if $\rho$ is an identity;
  we write $\strongDiag(\cat{X})$ and $\strictDiag(\cat{X})$ for the wide subcategories consisting of pseudo and strict morphisms, respectively.
\end{definition}

\begin{figure}[ht!]
  \centering
  \[
    \begin{tikzcd}[column sep=scriptsize]
      % https://q.uiver.app/?q=WzAsMyxbMCwwLCJcXG1hdGhiYntKfSJdLFsyLDAsIlxcbWF0aGJie0t9Il0sWzEsMSwiXFxtYXRoYmJ7WH0iXSxbMCwyLCJEIiwyXSxbMSwyLCJFIl0sWzAsMSwiUiJdLFszLDQsIlxccmhvIiwxLHsib2Zmc2V0IjotMiwic2hvcnRlbiI6eyJzb3VyY2UiOjIwLCJ0YXJnZXQiOjIwfX1dXQ==
      {\cat{J}} && {\cat{K}} \\
      & {\cat{X}}
      \arrow[""{name=0, anchor=center, inner sep=0}, "D"', from=1-1, to=2-2]
      \arrow[""{name=1, anchor=center, inner sep=0}, "E", from=1-3, to=2-2]
      \arrow["R", from=1-1, to=1-3]
      \arrow["\rho"{description}, shift left=2, shorten <=6pt, shorten >=6pt, Rightarrow, from=0, to=1]
    \end{tikzcd}
    \qquad\qquad\qquad
    \begin{tikzcd}[column sep=scriptsize]
      % https://q.uiver.app/?q=WzAsMyxbMCwwLCJcXG1hdGhiYntKfSJdLFsyLDAsIlxcbWF0aGJie0t9Il0sWzEsMSwiXFxtYXRoYmJ7WH0iXSxbMCwyLCJEIiwyXSxbMSwyLCJFIl0sWzEsMCwiUiIsMl0sWzMsNCwiXFxyaG8iLDEseyJvZmZzZXQiOi0yLCJzaG9ydGVuIjp7InNvdXJjZSI6MjAsInRhcmdldCI6MjB9fV1d
      {\cat{J}} && {\cat{K}} \\
      & {\cat{X}}
      \arrow[""{name=0, anchor=center, inner sep=0}, "D"', from=1-1, to=2-2]
      \arrow[""{name=1, anchor=center, inner sep=0}, "E", from=1-3, to=2-2]
      \arrow["R"', from=1-3, to=1-1]
      \arrow["\rho"{description}, shift left=2, shorten <=6pt, shorten >=6pt, Rightarrow, from=0, to=1]
    \end{tikzcd}
  \]
  \caption{\emph{Left:} a morphism in $\Diag(\cat{X})$. \emph{Right:} a morphism in $\DiagOp(\cat{X})$.}
  \label{fig:morphisms-in-diag-and-diag-op}
\end{figure}

\begin{remark}
  The construction of the diagram categories is functorial: given a functor $F\colon\cat{X}\to\cat{Y}$, we get a functor $\Diag(F)\colon\Diag(\cat{X})\to\Diag(\cat{Y})$ by post-composition with $F$, sending $(\cat{J},D)$ to $(\cat{J},F \circ D)$ and acting on morphisms $(R,\rho)$ by whiskering with $\rho.$ The same is true for $\DiagOp$.

  It is possible to give $\Diag$ and $\DiagOp$ a universal property as a kind of categorified slice category. Much as an ordinary slice category may be
  seen as the limit of a span in the 2-category $\Cat$,
  weighted by a certain span $\mathbb 1\to \mathbb 2\leftarrow \mathbb 1$ of categories, so can
  the diagram category be seen as the limit of a span
  in $\cat V$ weighted by the same span $\mathbb 1\to \mathbb 2\leftarrow \mathbb 1$ viewed in $\cat V$,
  where $\cat V$ is the category $\cat{OpLaxGray}$
  of $2$-categories enriched with the oplax Gray tensor product.
  We shall make no use of this fact, but point it out as
  an amusing reflection on the depth of vision of Eilenberg
  and Mac Lane's original paper, and perhaps as motivation
  for the work of this paper in getting \emph{out} of the
  relatively mysterious diagram category and into the more
  familiar slice $\Cat/\cat{X}.$
\end{remark}

Note that the wide subcategory $\strictDiag(\cat{X})$
consisting only of strict morphisms is precisely the slice
category $\Cat/\cat{X}$, and similarly for $\strictDiagOp(\cat
{X})$ and $(\Cat/\cat{X})^\op$. Thus $\strictDiagOp(\cat{X})
\cong\strictDiag(\cat{X})^\op$.
We can generalise this correspondence to the case of pseudo
morphisms using the following definition.

\begin{definition}
  The \emph{mate} of a morphism $(R,\rho)\colon D\to E$ in $\strongDiag(\cat{X})$ is the morphism $(R,\rho^{-1})\colon E\to D$ in $\strongDiagOp(\cat{X})$.
  We also refer to the inverse operation, which sends a morphism $(S,\sigma)$ in $\strongDiagOp(\cat{X})$ to the morphism $(S,\sigma^{-1})$ in $\strongDiag(\cat{X})$, as the \emph{mate}.
\end{definition}

It is immediately observed that the identity-on-objects functor $\strongDiag(\cat{X})\to\strongDiagOp(\cat{X})^\op$ that sends a morphism $(R,\rho)\colon D\to E$ to its mate $(R,\rho^{-1})\colon E\to D$ induces an isomorphism of categories.

This gives us one general motivation for studying diagram categories: they are extensions of the slice category $\Cat/\cat{X}$ that can describe more sensitive 2-categorical information. Objects of the diagram categories are also useful
as presentations for discrete (op)fibrations, as we shall see.
Furthermore, diagram categories are fundamental even at the 1-categorical level, as shown by the following lemma.

\begin{lemma}
\label{lemma:lim-gives-functor-diagopx-to-x}
  Let $\cat{X}$ be a complete category. Then giving a choice of limit for every diagram in $\cat{X}$ it is equivalent to giving a functor
  \[
    \lim\colon\DiagOp(\cat{X})\to\cat{X}.
  \]
  The dual statement for cocomplete categories also holds, giving a functor $\colim\colon\Diag(\cat{X})\to\cat{X}$.
\end{lemma}

\begin{proof}
  An early version of this statement appears in Eilenberg and Mac Lane's original paper on category theory \cite[§23]{eilenberg1945}; a proof can also be found in \cite[Proposition~4.2]{patterson2022}.
\end{proof}

The interest in $\Diag(\cat{X})$ and $\DiagOp(\cat{X})$ in \cite{patterson2022} was motivated by how diagrams can be used to present systems of equations, and how lifts of diagrams then describe solutions to these systems.
We will briefly recall the relevant definitions, but we refer the interested reader to \emph{loc.\@cit.\@} for a more detailed discussion.

\begin{definition}
\label{definition:lift-of-diagram}
  Given a diagram $D\colon\cat{J}\to\cat{X}$ and a functor $\pi\colon\cat{E}\to\cat{X}$, a \emph{lift of $D$ along $\pi$} is a functor $\overline{D}\colon\cat{J}\to\cat{E}$ making the triangle
  \[
    % https://q.uiver.app/#q=WzAsMyxbMSwwLCJcXGNhdHtFfSJdLFsxLDEsIlxcY2F0e1h9Il0sWzAsMSwiXFxjYXR7Sn0iXSxbMiwxLCJEIiwyXSxbMCwxLCJcXHBpIl0sWzIsMCwiXFxvdmVybGluZXtEfSIsMCx7InN0eWxlIjp7ImJvZHkiOnsibmFtZSI6ImRhc2hlZCJ9fX1dXQ==
    \begin{tikzcd}
      & {\cat{E}} \\
      {\cat{J}} & {\cat{X}}
      \arrow["D"', from=2-1, to=2-2]
      \arrow["\pi", from=1-2, to=2-2]
      \arrow["{\overline{D}}", dashed, from=2-1, to=1-2]
    \end{tikzcd}
  \]
  commutes.
  We write $\Lift{D}{\pi}$ for the set of lifts of $D$ along $\pi$.
\end{definition}

There is indeed a natural notion of \emph{category} of lifts along $\pi$,
where morphisms lie over $\id_D$, but for our purposes $\pi$ will always
be a discrete opfibration and thus the category of lifts will be discrete.

While not all of the following definitions are standard, they are consistent
with the usual ones: a functor is a discrete opfibration (in the usual sense)
precisely when it is a discrete opfibration at every morphism in $\cat{X}$.

\begin{definition}
\label{definition:discrete-opfibration}
  Let $\pi\colon\cat{E}\to\cat{X}$ be a functor.
  For any object $x\in\cat{X}$, we write $\cat{E}_x$ for the \emph{fibre of $\cat{E}$ over $x$}, defined as the set of all objects $\overline{x}\in\cat{E}$ such that $\pi(\overline{x})=x$.

  We say that $\pi$ is a \emph{discrete opfibration at a morphism $f\colon x \to y$} in $\cat{X}$ if, for every $\overline{x}\in\cat{E}_x$, there exists a unique object $\overline{y}\in\cat{E}$ and a unique morphism $\overline{f}\colon\overline{x}\to\overline{y}$ such that $\pi(\overline{f})=f$;
  we refer to $\overline{f}$ as the (unique) \emph{lift} of the pair $(\overline{x},f)$.

  Similarly, $\pi$ is said to be a \emph{discrete fibration at $f$} if $\pi^\op\colon\cat{E}^\op\to\cat{X}^\op$ is a discrete opfibration at $f$.
  If $\pi$ is both a discrete fibration at $f$ and a discrete opfibration at $f$, then we say that it is a \emph{discrete bifibration at $f$}.
\end{definition}

\begin{figure}[ht!]
  \centering
  \[
    \begin{tikzcd}
      % https://q.uiver.app/#q=WzAsNCxbMCwwLCJcXG92ZXJsaW5le3h9Il0sWzIsMCwiXFxidWxsZXQiXSxbMCwxLCJ4Il0sWzIsMSwieSJdLFsyLDMsImYiLDJdLFswLDEsIlxcb3ZlcmxpbmV7Zn0iLDAseyJzdHlsZSI6eyJib2R5Ijp7Im5hbWUiOiJkYXNoZWQifX19XSxbNSw0LCJcXHBpIiwwLHsic2hvcnRlbiI6eyJzb3VyY2UiOjIwLCJ0YXJnZXQiOjIwfSwibGV2ZWwiOjEsInN0eWxlIjp7InRhaWwiOnsibmFtZSI6Im1hcHMgdG8ifX19XV0=
      {\overline{x}} && \bullet \\
      x && y
      \arrow[""{name=0, anchor=center, inner sep=0}, "f"', from=2-1, to=2-3]
      \arrow[""{name=1, anchor=center, inner sep=0}, "{\overline{f}}", dashed, from=1-1, to=1-3]
      \arrow["\pi", shorten <=4pt, shorten >=4pt, maps to, from=1, to=0]
    \end{tikzcd}
  \]
  \caption{A discrete opfibration $\pi$ at $f\colon x\to y$. As the picture suggests, the codomain of the lift $\overline{f}$ is not given, but is instead part of the existence and uniqueness statement.}
\end{figure}

\begin{remark}
\label{remark:disc-opfib-functorial-lifting}
  An important consequence of the uniqueness of lifts is that discrete opfibrations give \emph{functorial} lifts: given morphisms $f\colon x\to y$ and $g\colon y\to z$ in $\cat{X}$ and an object $\overline{x}$ in $\cat{E}_x$, we can lift $(\overline{x},f)$ to get $\overline{f}\colon\overline{x}\to\overline{y}$, and then lift $(\overline{y},g)$ to get $\overline{g}\colon\overline{y}\to\overline{z}$; by uniqueness of lifts, it must be the case that $\overline{g}\circ\overline{f}=\overline{g\circ f}$, since the former also defines a lift of $g\circ f\colon x\to z$. This is at the heart of the proof that the Grothendieck construction gives as equivalence between the category of copresheaves on $\cat{X}$ and the category of discrete opfibrations over $\cat{X}.$
\end{remark}

\begin{remark}
  We have so far used the word ``lift'' to mean two seemingly different things: in \cref{definition:lift-of-diagram} it means a functor making a triangle commute, and in \cref{definition:discrete-opfibration} it means an object or morphism in the fibre of a functor.
  But the former is actually an example of the latter, since finding a lift of $D\colon\cat{J}\to\cat{X}$ along $\pi\colon\cat{E}\to\cat{X}$ is exactly finding an object in the fibre of $(\pi\circ-)\colon[\cat{J},\cat{E}]\to[\cat{J},\cat{X}]$ over $D$.
  With this apology made, we shall feel free to use the same overline notation for lifts in both senses.
\end{remark}

The perspective of diagrams presenting systems of equations leads to a question, however, that remains largely invisible in standard algebraic notation for equations: there can exist diagrams that are not isomorphic but which nonetheless present systems of equations that have ``the same'' sets of solutions.
The authors of \cite{patterson2022} propose a definition of a \emph{weak equivalence} of diagrams, which aims to rectify this situation by defining two diagrams to be weakly equivalent exactly when there is a pseudo map between them
inducing a bijection between sets of lifts against all discrete opfibrations. As discussed in the introduction, these lifting sets are interpreted as the sets
of solutions to the systems of equations presented by the two
diagrams.

In this paper we somewhat weaken this definition,
allowed for weak equivalences that aren't necessarily pseudo,
and avoiding dependence on the mate construction in the
definition. In the below definition we consider only weak
equivalences in $\DiagOp(\cat{X})$, as it is this category
whose morphisms induce functions between lifting sets against
discrete opfibrations.

\begin{definition}
\label{definition:equivalence-of-systems}
  We say that a morphism $(R,\rho)\colon D\to E$ in $\DiagOp(\cat{X})$ is a \emph{weak equivalence} if, for every discrete opfibration $\pi\colon\cat{E}\to\cat{X}$, the map $R_*\colon\Lift{D}{\pi}\to\Lift{E}{\pi}$ induced by $(R,\rho)$ (see below) is a bijection.
\end{definition}

The map $R_*$ is given by applying the discrete opfibration
property of $\pi$ to find a lift of $\rho$ whose codomain is
the image of a lift of $D$, as we now explain in more detail,
following \cite[Theorem 5.2]{patterson2022}. A slightly
different proof is given at a higher level of abstraction in \cref{section:2-categorical-story}.

For any $k\in\cat{K}$, the natural transformation $\rho$ gives us a morphism
\[
  \rho_k\colon DR(k)\to E(k)
\]
in $\cat{X}$.
If we have some lift $\overline{D}\colon\cat{J}\to\cat{E}$ of $D\colon\cat{J}\to\cat{X}$ along $\pi$, then we have the necessary data to apply the discrete opfibration property: there exists a unique lift $\overline{\rho_k}$ of $\rho_k$, which has a codomain we shall name $\overline{E(k)}$, as seen here:
\[
  \begin{tikzcd}
    % https://q.uiver.app/#q=WzAsNCxbMCwwLCJcXG92ZXJsaW5le0R9UihrKSJdLFsyLDAsIlxcb3ZlcmxpbmV7RShrKX0iXSxbMCwxLCJEUihrKSJdLFsyLDEsIkUoaykiXSxbMiwzLCJcXHJob19rIiwyXSxbMCwxLCJcXG92ZXJsaW5le1xccmhvX2t9IiwwLHsic3R5bGUiOnsiYm9keSI6eyJuYW1lIjoiZGFzaGVkIn19fV0sWzUsNCwiXFxwaSIsMCx7InNob3J0ZW4iOnsic291cmNlIjoyMCwidGFyZ2V0IjoyMH0sImxldmVsIjoxLCJzdHlsZSI6eyJ0YWlsIjp7Im5hbWUiOiJtYXBzIHRvIn19fV1d
    {\overline{D} R(k)} && {\overline{E(k)}} \\
    {DR(k)} && {E(k)}
    \arrow[""{name=0, anchor=center, inner sep=0}, "{\rho_k}"', from=2-1, to=2-3]
    \arrow[""{name=1, anchor=center, inner sep=0}, "{\overline{\rho_k}}", dashed, from=1-1, to=1-3]
    \arrow["\pi", shorten <=4pt, shorten >=4pt, maps to, from=1, to=0]
  \end{tikzcd}
\]
With that done, we define the functor $\overline{E}\colon\cat{K}\to\cat{E}$ on objects by $\overline{E}(k)=\overline{E(k)}$.
Then, given a morphism $g\colon k\to k'$ in $\cat{K}$, we can apply the discrete opfibration property again: there exists a unique lift $\overline{E(g)}$ of domain $\bar E(k).$ We call its codomain $\widetilde{E(k')}:$
\[
  \begin{tikzcd}
    % https://q.uiver.app/#q=WzAsNCxbMCwwLCJcXG92ZXJsaW5le0V9KGspIl0sWzIsMCwiXFx3aWRldGlsZGV7RShrJyl9Il0sWzAsMSwiRShrKSJdLFsyLDEsIkUoaycpIl0sWzIsMywiRShnKSIsMl0sWzAsMSwiXFxvdmVybGluZXtFKGcpfSIsMCx7InN0eWxlIjp7ImJvZHkiOnsibmFtZSI6ImRhc2hlZCJ9fX1dLFs1LDQsIlxccGkiLDAseyJzaG9ydGVuIjp7InNvdXJjZSI6MjAsInRhcmdldCI6MjB9LCJsZXZlbCI6MSwic3R5bGUiOnsidGFpbCI6eyJuYW1lIjoibWFwcyB0byJ9fX1dXQ==
    {\overline{E}(k)} && {\widetilde{E(k')}} \\
    {E(k)} && {E(k')}
    \arrow[""{name=0, anchor=center, inner sep=0}, "{E(g)}"', from=2-1, to=2-3]
    \arrow[""{name=1, anchor=center, inner sep=0}, "{\overline{E(g)}}", dashed, from=1-1, to=1-3]
    \arrow["\pi", shorten <=4pt, shorten >=4pt, maps to, from=1, to=0]
  \end{tikzcd}
\]
In order to define $\overline{E}(g)=\overline{E(g)}$, we need to know that $\widetilde{E(k')}=\bar E(k').$
This follows directly from discrete opfibrations having functorial lifting (\cref{remark:disc-opfib-functorial-lifting}), but for clarity we explain the details in this case.
We already know that $\overline{\rho_k}$ is the unique lift of $\rho_k$ with domain $\bar D R(k).$ Similarly, $\overline{E(g)}$ is by definition the unique lift of $E(g)$ with domain $\bar E(k).$

Together, this implies that $\overline{E(g)}\circ\overline{\rho_k}$,  which has codomain $\widetilde{E(k')}$, is the unique lift of $E(g)\circ\rho_k$ with domain $\overline{D} R(k).$
But the naturality of $\rho$ says that $E(g)\circ\rho_k=\rho_{k'}\circ DR(g)$, and by a similar argument as above, the unique lift of $\rho_{k'}\circ DR(g)$ with domain $\overline{D} R(k)$ has codomain $\overline{E}(k).$ Thus  $\widetilde{E(k')}=\overline{E(k')}$, as desired.

\begin{remark}
  Note that this argument breaks down if $R$ goes in the other direction: morphisms in $\Diag(\cat{X})$ do not necessarily induces maps of lifts against discrete opfibrations.
  For this reason, here and in \cite{patterson2022} we have been primarily interested in the category $\DiagOp(\cat{X})$.
\end{remark}

In the above construction of $R_*\colon\Lift{D}{\pi}\to\Lift{E}{\pi}$, we have actually also constructed a natural transformation $\overline{\rho}\colon\overline{D}R\Rightarrow\overline{E}$ by setting $\overline{\rho}_k=\overline{\rho_k}$, which gives us a morphism $(R,\overline{\rho})\colon\overline{D}\to\overline{E}$ in $\DiagOp(\cat{E})$.
In fact, by the uniqueness of the maps constructed throughout, we have actually established the following result.

\begin{lemma}[{\cite[Theorem 5.2]{patterson2022}}]
\label{lemma:diagop-disc-opfib-is-disc-opfib}
  The functor $\DiagOp$ preserves discrete opfibrations. That is, let $\pi\colon\cat{E}\to\cat{X}$ be a discrete opfibration.
  Then the post-composition functor
  \[
    \DiagOp(\pi)\colon\DiagOp(\cat{E})\to\DiagOp(\cat{X})
  \]
  is also a discrete opfibration.
\end{lemma}

We can even express the definition of a weak equivalence in $\DiagOp(\cat{X})$ purely in terms of properties of the functor $\DiagOp(\pi)$, avoiding explicit mention of lifts of diagrams altogether.
To do so, we need to give an alternative formulation of discrete opfibrations.

\begin{definition}
\label{definition:discrete-opfibration-induces-map-on-fibres}
  Let $\pi\colon\cat{E}\to\cat{X}$ be a discrete opfibration at some morphism $f\colon x\to y$ in $\cat{X}$.
  The discrete opfibration property of $\pi$ induces a function between fibres which we denote by
  \[
    \begin{aligned}
      f_*\colon\cat{E}_x
      &\to\cat{E}_y
    \\\overline{x}
      &\mapsto\overline{y}
    \end{aligned}
  \]
  that sends $\overline{x}$ to the codomain $\bar y$ of the unique lift of $f$ with domain $\bar x.$

  Dually, if $\pi$ is a discrete fibration at $f$, then it induces a function
  \[
    \begin{aligned}
      f^*\colon\cat{E}_y
      &\to\cat{E}_x
    \\\overline{y}
      &\mapsto\overline{x}
    \end{aligned}
  \]
  that sends $\overline{y}$ to the domain of the unique lift $\overline{f}$ of $f$ with codomain $\bar y.$
\end{definition}

\begin{lemma}
\label{lemma:discrete-bifibrations-and-bijections}
  Let $\pi\colon\cat{E}\to\cat{X}$ be a discrete opfibration at some morphism $f\colon x\to y$ in $\cat{X}$.
  Then $\pi$ is a discrete fibration (and thus a discrete bifibration) at $f$ if and only if the induced map $f_*\colon\cat{E}_x\to\cat{E}_y$ is a bijection.
\end{lemma}

\begin{proof}
  If $\pi$ is also a discrete fibration at $f$, then we have the induced map $f^*\colon\cat{E}_y\to\cat{E}_x$ which sends $\overline{y}$ to the domain of the corresponding unique lift of $f$.
  But the uniqueness of lifts implies that $f^*f_*=\id_{\cat{E}_x}$ and $f_*f^*=\id_{\cat{E}_y}$ and so $f_*$ is a bijection.

  Conversely, if $f_*$ is a bijection, then define $f^*\colon\cat{E}_y\to\cat{E}_x$ by $f^*(\overline{y})=f_*^{-1}(\overline{y})$.
  Then $f^*$ gives us the necessary domain for the discrete fibration property: let $\overline{y}\in\cat{E}_y$ and set $\overline{x}=f^*(\overline{y})$; since $f_*(\overline{x})=\overline{y}$ by construction, the discrete opfibration property gives a unique lift $\overline{f}\colon\overline{x}\to\overline{y}$.
\end{proof}

\begin{corollary}
\label{corollary:equivalence-iff-disc-opfib-to-disc-opfib}
  A morphism $(R,\rho)$ in $\DiagOp(\cat{X})$ is a weak equivalence if and only if, for every discrete opfibration $\pi\colon\cat{E}\to\cat{X}$, the functor $\DiagOp(\pi)\colon\DiagOp(\cat{E})\to\DiagOp(\cat{X})$ is a discrete fibration (and thus a discrete bifibration, by \cref{lemma:diagop-disc-opfib-is-disc-opfib}) at $(R,\rho)$.
\end{corollary}

\begin{proof}
  This follows simply by unravelling the definitions.
  First of all, to combine the notation of \cref{definition:discrete-opfibration-induces-map-on-fibres,definition:equivalence-of-systems},
  notice that the objects of $\DiagOp(\cat{E})_D$ are exactly the elements of $\Lift{D}{\pi}$, and a morphism $(R,\rho)\colon D\to E$ in $\DiagOp(\cat{X})$ induces the function on fibres
  \[
    \begin{aligned}
      R_*\colon \DiagOp(\cat{E})_D
      &\to \DiagOp(\cat{E})_E
    \\\overline{D}
      &\mapsto \cod(R,\overline{\rho})
    \end{aligned}
  \]
  where $\overline{\rho}$ is constructed as in the proof of \cref{lemma:diagop-disc-opfib-is-disc-opfib}, and $\cod$ denotes the codomain.
  Note that this map $R_*$ is exactly the map $R_*$ in \cref{definition:equivalence-of-systems}.

  Then, by \cref{lemma:discrete-bifibrations-and-bijections}, we know that $\pi$ is a discrete fibration at $(R,\rho)$ if and only if this map $R_*$ is bijective, but this is exactly the the definition of a weak equivalence.
\end{proof}

\begin{example}
\label{example:representable-dopfs-do-no-suffice}
	Every discrete opfibration over $\cat{X}$ can be constructed out of the ``representable'' discrete opfibrations
	$\pi_x:x/\cat{X}\to \cat{X}$, where $x$ is an object of $\cat{X}$, via colimits. As described in the introduction,
	lifting a diagram $D$ along $\pi_x$ can be interpreted as giving an $x$-parameterized family of solutions to the system
	of equations presented by $D$.
	Since this is the most familiar notion of solution to a system of equations, it is natural
	to ask whether our notion of weak equivalence depends on the use of arbitrary discrete opfibrations, and in particular
	on whether there are non-weakly-equivalent diagrams admitting a map that does induce a bijection on lifts against
	representable discrete opfibrations.

	Indeed, this can happen, intuitively because we are mapping \emph{into} discrete opfibrations when constructing lifts but
	discrete opfibrations are built from representables via \emph{colimits}, not via limits. For a minimal example, let $\cat{X}$
	be the arrow category $\mathbb 2= (0\to 1).$ Then representable discrete opfibrations are given by the identity and by 
  the inclusion $1\colon\mathbb{1}\to \cat X.$ Every morphism of diagrams induces a bijection on lifts 
  against the identity; for the nontrivial representable discrete opfibration, a lift of $D\colon J\to \cat X$ along
  $1$ exists if and only if $D$ is constant at $1,$ in which case this lift is unique. Thus morphisms of 
  diagrams in $\cat X$ which are either both constant at $1$ or neither constant at $1$ induce bijections on
  lifts against representable discrete opfibrations. However, there is a unique diagram morphism 
  $1\to 1+1$ determined by the unique functor $\mathbb{1}+\mathbb{1}\to \mathbb{1},$ and it certainly does
  not induce a bijection on lifts against $1+1$ itself, since $1+1$ has four lifts along itself but
  $1$ has only two lifts along $1+1.$
\end{example}

We now turn to the other diagram category, $\Diag(\cat{X})$.
As mentioned above, it is not necessarily the case that a morphism in $\Diag(\cat{X})$ induces a map on lifts along a discrete opfibration, so we give a name to those that do.

\begin{definition}
\label{definition:pushforward-property}
  We say that a morphism $(R,\rho)$ in $\Diag(\cat{X})$ is a \emph{weak equivalence} if, for any discrete opfibration $\pi\colon\cat{E}\to\cat{X}$, the functor $\Diag(\pi)\colon\Diag(\cat{E})\to\Diag(\cat{X})$ is a discrete opfibration at $(R,\rho)$.
\end{definition}

The dual version of \cref{lemma:diagop-disc-opfib-is-disc-opfib} is true in a limited sense: if we consider only pseudo morphisms, then $\Diag$ sends discrete opfibrations to discrete fibrations.
More precisely, we have the following lemma.

\begin{lemma}
\label{lemma:strongdiag-disc-opfib-is-disc-fib}
  Let $\pi\colon\cat{E}\to\cat{X}$ be a discrete opfibration.
  Then
  \[
    \strongDiag(\pi)\colon\strongDiag(\cat{E})\to\strongDiag(\cat{X})
  \]
  is a discrete fibration.
\end{lemma}

\begin{proof}
  The idea of the proof is rather simple --- take the mate of a morphism in $\strongDiag(\cat{X})$ and then apply the proof of \cref{lemma:diagop-disc-opfib-is-disc-opfib}.

  We need to show that, given some (pseudo) morphism $(S,\sigma)\colon(D,\cat{J})\to(E,\cat{K})$ in $\strongDiag(\cat{X})$ and a lift $\overline{E}$ of $E$ along $\pi$, there exists a unique lift $(S,\overline{\sigma})\colon(\overline{D},\cat{J})\to(\overline{E},\cat{K})$ in $\strongDiag(\cat{E})$ of $(S,\sigma)$.
  For any $j\in\cat{J}$, we have the morphism
  \[
    \sigma_j\colon D(j)\to ES(j)
  \]
  in $\cat{X}$.
  By hypothesis, $\sigma$ is a natural isomorphism, which means that $\sigma_j$ is an isomorphism, and thus we have an inverse morphism
  \[
    \sigma_j^{-1}\colon ES(j)\to D(j).
  \]
  But now we are in exactly the same situation as in the proof of \cref{lemma:diagop-disc-opfib-is-disc-opfib}, just with the roles of $D$ and $E$ swapped, and so we can apply exactly the same argument to obtain
  \[
    (S,\overline{\sigma^{-1}})\colon(\overline{E},\cat{K})\to(\overline{D},\cat{J})
  \]
  in $\strongDiagOp(\cat{E})$.
  Since $\overline{\sigma^{-1}}$ thus defined is a natural isomorphism, we can again take its inverse in order to obtain
  \[
    (S,\overline{\sigma})\colon(\overline{D},\cat{J})\to(\overline{E},\cat{K})
  \]
  in $\strongDiag(\cat{E})$, as desired.
\end{proof}

We now know that, if $\pi\colon\cat{E}\to\cat{X}$ is a discrete opfibration, then $\DiagOp(\pi)$ (and thus $\strongDiagOp(\pi)$, as we will explain in the proof of \cref{corollary:pushforward-iff-equivalence}) is also a discrete opfibration, and $\strongDiag(\pi)$ is a discrete fibration.
It is natural to then ask if $\strongDiagOp(\pi)$ being a discrete fibration is equivalent to $\strongDiag(\pi)$ being a discrete opfibration: do weak equivalences in $\strongDiagOp(\cat{X})$ and weak equivalences in $\strongDiag(\cat{X})$ coincide, up to the mate isomorphism?
The following lemma and its corollary provides an affirmative answer to this question.

\begin{lemma}
\label{lemma:strongdiagop-disc-fib-iff-strongdiag-disc-opfib}
  Let $\pi\colon\cat{E}\to\cat{X}$ be a discrete opfibration, and let $(R,\rho)\colon D\to E$ be a morphism in $\strongDiagOp(\cat{X})$.
  Then $\strongDiagOp(\pi)$ is a discrete fibration at $(R,\rho)$ if and only if $\strongDiag(\pi)$ is a discrete opfibration at $(R,\rho^{-1})$.
\end{lemma}

\begin{proof}
  By \cref{lemma:discrete-bifibrations-and-bijections}, we know that $\strongDiagOp(\pi)$ is a discrete fibration at $(R,\rho)$ if and only if $(R,\rho)_*\colon\Lift{\pi}{D}\to \Lift{\pi}{E}$ is a bijection.
  The dual statement of the same lemma tells us that $\strongDiag(\pi)$ is a discrete opfibration at $(R,\rho^{-1})$ if and only if $(R,\rho^{-1})^*\colon\Lift{\pi}{D}\to \Lift{\pi}{E}$ is a bijection. Thus we need only show that $(R,\rho)_*$ and $(R,\rho^{-1})^*$ are the same function.

  By definition, given $\bar D\in \Lift{\pi}{D}$, the image $(R,\rho)_*(\bar D)$ is the codomain of the unique lift of $(R,\rho)$ with domain $\bar D$, while $(R,\rho^{-1})^*(\bar D)$ is the domain of the unique lift of $(R,\rho^{-1})$ with codomain $\bar D.$
  Since the discrete opfibration $\pi$ is conservative, in a lift $(R,\bar\rho)\colon\bar D\to \bar E$ of $(R,\rho)$ we can be sure that
  $\bar\rho$ is invertible. Thus the mate $(R,\bar\rho^{-1})\colon\bar E\to \bar D$ of the unique such lift must also be the
  unique lift determining $(R,\rho^{-1})^*.$ Thus $(R,\rho)_*$ and $(R,\rho^{-1})^*$ are the same function, as desired.
\end{proof}

\begin{corollary}
\label{corollary:pushforward-iff-equivalence}
  A morphism $(R,\rho)\colon D\to E$ in $\strongDiagOp(\cat{X})$ is a weak equivalence in the sense of \cref{definition:equivalence-of-systems} if and only if its mate $(R,\rho^{-1})\colon E\to D$ in $\strongDiag(\cat{X})$ is a weak equivalence in the sense of \cref{definition:pushforward-property}.
\end{corollary}

\begin{proof}
  We know from \cref{corollary:equivalence-iff-disc-opfib-to-disc-opfib} that $(R,\rho)$ is a weak equivalence if and only if, for any discrete opfibration $\pi\colon\cat{E}\to\cat{X}$, the discrete opfibration $\DiagOp(\pi)\colon\DiagOp(\cat{E})\to\DiagOp(\cat{X})$ is also a discrete fibration at $(R,\rho)$.
  But any lift $(R,\overline{\rho})$ of $(R,\rho)$ along $\pi$ will be a pseudo morphism in $\DiagOp(\cat{E})$, since discrete opfibrations lift isomorphisms to isomorphisms (a consequence of functorial lifting, as in \cref{remark:disc-opfib-functorial-lifting}) and so $\overline{\rho}$ will be a natural isomorphism since $\rho$ is.
  This means that $\DiagOp(\pi)$ is a discrete fibration at $(R,\rho)$ if and only if $\strongDiagOp(\pi)$ is a discrete fibration at $(R,\rho)$.
  But then \cref{lemma:strongdiagop-disc-fib-iff-strongdiag-disc-opfib} tells us that this is equivalent to $\strongDiag(\pi)$ being a discrete opfibration at $(R,\rho^{-1})$, which is exactly the definition of $(R,\rho^{-1})$ being a weak equivalence.
\end{proof}

It turns out that there is a readily verifiable sufficient condition on the functor $R$ that ensures that a morphism $(R,\rho)$ in $\Diag$ is a weak equivalence.
That condition is \emph{initiality}, which we now review.

\section{Initial and relatively initial functors}
\label{section:initial-and-relatively-functors}

In this section we will recall the definition of an initial functor in the context of orthogonal factorisation systems, and then explain how this relates to the weak equivalences from the previous section.
This motivates a generalisation of initiality, leading to the definition of relatively initial morphisms.
Note that throughout this section we will be working exclusively with $\Diag$, not $\DiagOp$.

Informally, a \emph{zigzag} between two objects $x$ and $x'$, denoted $x\zigzag x'$, is a finite sequence of morphisms $x=x_0\leftrightarrow x_1\leftrightarrow x_2\leftrightarrow\ldots\leftrightarrow x_{n+1}=x'$ in which each arrow is allowed to point in either direction (``$\leftrightarrow\ \in\{\from,\to\}$''). More formally:

\begin{definition}
  A \emph{zigzag} between objects $x$ and $x'$ in a category $\cat{X}$ is a finite sequence of spans $\{x_i\from y_i\to x_{i+1}\}_{i=0,\ldots,n}$ in $\cat{X}$ such that $x_0=x$ and $x_{n+1}=x'$.
\end{definition}

We can recover the informal definition of a zigzag from the formal one by taking some of the legs of the spans $x_i\from y_i\to x_{i+1}$ to be identities, so that the directions of the arrows don't appear to strictly alternate.

\begin{definition}
  A functor $R\colon\cat{J}\to\cat{K}$ is \emph{initial} if, for every object $k\in\cat{K}$, the comma category $R/k$ is \emph{non-empty} and \emph{connected}, that is, for every $k\in\cat{K}$,
  \begin{itemize}
    \item \emph{(non-empty)} there exists at least one object $j\in\cat{J}$ along with a morphism $f\colon Rj\to k$ in $\cat{K}$; and
    \item \emph{(connected)} given any two objects $(j,f)$ and $(j',f')$ in $R/k$, there exists a zigzag $(j,f)\zigzag(j',f')$ in $R/k$.
  \end{itemize}
  More abstractly, a category $\cat{X}$ is \emph{connected} if the right adjoint $\pi_0\colon\Cat \to \Set$ to the chaotic/codiscrete category functor $\Set\to\Cat$ sends $\cat{X}$ to a point.
\end{definition}

A key way of thinking about initial functors is through
following characterization, that restricting along an initial functor does not modify a diagram's limit.

\begin{lemma}
  Let $R\colon\cat{J}\to\cat{K}$ be a functor.
  Then $R$ is initial if and only if, for any category $\cat{X}$ and any functor $D\colon\cat{J}\to\cat{X}$, the canonical morphism
  \[
    \lim_{\cat{K}}D \to \lim_{\cat{J}}(DR)
  \]
  is an isomorphism whenever the limits involved exist.
  Furthermore, it is sufficient to establish this in the case that $\cat{X}=\Set$.
\end{lemma}

\begin{proof}
  See, for example, \cite[Lemma~8.3.4]{riehl2014}.
\end{proof}

Two particularly useful sufficient conditions for initiality are the following.

\begin{lemma}
\label{lemma:left-adjoints-are-initial}
  Let $L\colon\cat{J}\to\cat{K}$ be a left adjoint.
  Then $L$ is initial.
\end{lemma}

\begin{proof}
  One way of defining what it means for $L$ to be a
  left adjoint is that the comma category $L/k$ has
  a terminal object $\varepsilon_k:LRk\to k$ for every
  $k.$ Any category with a terminal object is \emph{a fortiori} connected.
\end{proof}

\begin{lemma}
\label{lemma:full-and-essentially-surjective-is-initial}
  Let $R\colon\cat{J}\to\cat{K}$ be a full and essentially surjective functor.
  Then $R$ is initial.
\end{lemma}

\begin{proof}
  This is \cite[Lemma~9.8]{patterson2022}.
\end{proof}

The fact that a morphism $(R,\rho)$ in $\Diag(\cat{X})$ is a weak equivalence if $R$ is initial is related to how initial functors and discrete opfibrations interact with each other.
We can make this more precise by recalling the so-called
comprehensive factorization system.

\begin{definition}
\label{definition:orthogonal-factorisation-system}
  Let $\cat{L}$ and $\cat{R}$ be classes of morphisms in a category $\cat{C}$.
  We say that $(\cat{L},\cat{R})$ is an \emph{orthogonal factorisation system} if
  \begin{itemize}[noitemsep]
    \item every morphism $f\colon c\to d$ in $\cat{C}$ factors as $f = r \circ l$ for some morphism $l \in \cat{L}$ and $r \in \cat{R}$, and the factorisation is unique up to unique isomorphism:
    \[
      % https://q.uiver.app/#q=WzAsNCxbMCwxLCJjIl0sWzIsMSwiZCJdLFsxLDAsImUiXSxbMSwyLCJlJyJdLFswLDIsImwiXSxbMiwxLCJyIl0sWzAsMywibCciLDJdLFsyLDMsIlxcZXhpc3RzICEiLDIseyJzdHlsZSI6eyJib2R5Ijp7Im5hbWUiOiJkYXNoZWQifX19XSxbMywxLCJyJyIsMl1d
      \begin{tikzcd}
        & e \\
        c && d \\
        & {e'}
        \arrow["l", from=2-1, to=1-2]
        \arrow["r", from=1-2, to=2-3]
        \arrow["{l'}"', from=2-1, to=3-2]
        \arrow["{\exists !}"', dashed, from=1-2, to=3-2]
        \arrow["{r'}"', from=3-2, to=2-3]
      \end{tikzcd}
      \qquad\text{whenever}\qquad r \circ l = f = r' \circ l';
    \]
    \item $\cat{L}$ and $\cat{R}$ contains all isomorphisms;
    \item $\cat{L}$ and $\cat{R}$ are closed under composition.
  \end{itemize}
  We refer to $\cat{L}$ as the \emph{left part} of the factorisation system, and to $\cat{R}$ as the \emph{right part}.
\end{definition}

\begin{remark}
  \label{remark:lifting-properties-for-ofs}
  A well-known equivalent definition says that classes of morphisms $(\cat L,\cat R)$ in a category $\cat{C}$ form an orthogonal factorization
  system if:
  \begin{itemize}
    \item Every morphism in $\cat C$ factors in \emph{some} way as an $\cat L$-map followed by an $\cat R$-map; and
    \item $\cat L$ and $\cat R$ are \emph{orthogonal} in the sense that,
  whenever a solid square as below is given with $\ell\in \cat L$ and $r\in \cat R$, there
  exists a unique diagonal morphism indicated by the dashed arrow making
  both triangles commute.
  \end{itemize}
  \begin{equation*}
    % https://q.uiver.app/#q=WzAsNCxbMCwwLCJjXzEiXSxbMCwxLCJkXzEiXSxbMSwwLCJjXzIiXSxbMSwxLCJkXzIiXSxbMCwxLCJcXGVsbCIsMl0sWzIsMywiciJdLFswLDJdLFsxLDNdLFsxLDIsIiIsMSx7InN0eWxlIjp7ImJvZHkiOnsibmFtZSI6ImRvdHRlZCJ9fX1dXQ==
    \begin{tikzcd}
      {c_1} & {c_2} \\
      {d_1} & {d_2}
      \arrow["\ell"', from=1-1, to=2-1]
      \arrow["r", from=1-2, to=2-2]
      \arrow[from=1-1, to=1-2]
      \arrow[from=2-1, to=2-2]
      \arrow[dotted, from=2-1, to=1-2]
    \end{tikzcd}
  \end{equation*}
  The proof is straightforward.
  Assuming the assumptions just listed, if a morphism $f\colon c_1\to d_2$ factors in two ways
  as an $\ell$ followed by an $r$, then that produces an orthogonality situation in two different
  ways, producing arrows $d_1 \leftrightarrows c_2$ which are mutually inverse isomorphisms by
  the usual arguments involving universal properties.
  \begin{equation*}
    % https://q.uiver.app/#q=WzAsNCxbMCwwLCJjXzEiXSxbMCwxLCJkXzEiXSxbMSwwLCJjXzIiXSxbMSwxLCJkXzIiXSxbMCwxLCJcXGVsbF8xIiwyXSxbMiwzLCJyXzEiXSxbMCwyLCJcXGVsbF8yIl0sWzEsMywicl8yIiwyXSxbMSwyLCIiLDEseyJzdHlsZSI6eyJ0YWlsIjp7Im5hbWUiOiJhcnJvd2hlYWQifSwiYm9keSI6eyJuYW1lIjoiZG90dGVkIn19fV1d
    \begin{tikzcd}
      {c_1} & {c_2} \\
      {d_1} & {d_2}
      \arrow["{\ell_1}"', from=1-1, to=2-1]
      \arrow["{r_1}", from=1-2, to=2-2]
      \arrow["{\ell_2}", from=1-1, to=1-2]
      \arrow["{r_2}"', from=2-1, to=2-2]
      \arrow[dotted, tail reversed, from=2-1, to=1-2]
    \end{tikzcd}
  \end{equation*}
  Conversely, assuming the conditions in \cref{definition:orthogonal-factorisation-system},
  when given an orthogonality situation as above, one factors the horizontal arrows
  to produce two factorizations of the path around the square and then uses essential
  uniqueness of the factorization to produce the desired diagonal filler.
  \end{remark}

The orthogonal factorisation system most relevant to us is the comprehensive factorisation \cite{street1973}, which will prove useful in the comparison between initial functors and weak equivalences in $\Cat/\cat{X}$ in \cref{section:localisations-of-categories-of-diagrams}.

\begin{lemma}
\label{lemma:comprehensive-factorisation-system}
  The category $\Cat$ admits an orthogonal factorisation system $(\cat L,\cat R)$, known as the \emph{comprehensive factorisation system}, where $\cat L$ is the class of initial functors and $\cat R$ is the class of discrete opfibrations.
\end{lemma}

\begin{proof}
  This is \cite[Theorems~3 and 4]{street1973}, but it is interesting to note that it also follows from the proof of \cref{lemma:relatively-initial-satisfies-pushforward} (see \cite[Theorem~9.12]{patterson2022}). The idea of the factorization is to use the Yoneda embedding to interpret
  $R$ as a diagram of copresheaves on $\cat K$ and take the colimit
  of this diagram, followed by the Grothendieck construction,
  to get the canonical discrete opfibration over $\cat K$
  through which $R$ factors.
\end{proof}

The following lemma provides the link between initial functors and weak equivalences in $\Diag$, keeping \cref{corollary:pushforward-iff-equivalence} in mind.

\begin{lemma}
\label{lemma:initial-satisfies-pushforward}
  Let $(R,\rho)$ be a morphism in $\Diag(\cat{X})$.
  If $R$ is an initial functor, then $(R,\rho)$ is a weak equivalence.
\end{lemma}

\begin{proof}
  This is a special case of
  \cref{lemma:relatively-initial-satisfies-pushforward} stated below.
\end{proof}

Although \cref{lemma:initial-satisfies-pushforward} tells us that initiality is a \emph{sufficient} condition for a morphism to be a weak equivalence, it is not \emph{necessary}, as demonstrated by \cref{example:relatively-initial-but-not-initial} below.
Indeed, this condition does not even make reference to the diagrams themselves, nor to the transformation $\rho$, so it is not surprising that it is too strong.

This defect is partially remedied by a more subtle version of initiality that does take into account both of these aspects.

\begin{definition}
\label{definition:relative-slice-category}
  Let $(R,\rho)\colon(\cat{J},D)\to(\cat{K},E)$ be a morphism in $\Diag(\cat{X})$, and let $k$ be an object of $\cat{K}$.
  Then we define the \emph{relative comma category} $(R,\rho)/k$ as follows:
  \begin{itemize}
    \item its objects are pairs $(j,f)$, where $j\in\cat{J}$ and $f\colon Rj\to k$ in $\cat{K}$;
    \item its morphisms $(j,f)\to(j',f')$ are morphisms $h\colon j\to j'$ in $\cat{J}$ such that the diagram
      \[
        \begin{tikzcd}[column sep=tiny]
          Dj
            \ar[rr,"Dh"]
            \ar[d,swap,"\rho_j"]
          &&Dj'
            \ar[d,"\rho_{j'}"]
        \\ERj
            \ar[dr,swap,near start,"Ef"]
          &&ERj'
            \ar[dl,near start,"Ef'"]
        \\&Ek
        \end{tikzcd}
      \]
      commutes.
  \end{itemize}
  We sometimes write the objects $(j,f)$ simply as $Rj\xrightarrow{f}k$.
\end{definition}

The relative comma category $(R,\rho)/k$ has the same objects as the comma category $R/k$, but has, in general, more morphisms.
Indeed, any morphism $(j,f)\to(j',f')$ in $R/k$, given by the data of a morphism $h\colon j\to j'$ in $\cat{J}$, defines a morphism $(j,f)\to(j',f')$ in $(R,\rho)/k$, since the naturality of $\rho\colon D\To ER$ implies that the diagram
\[
  \begin{tikzcd}[column sep=tiny]
    Dj
      \ar[rr,"Dh"]
      \ar[d,swap,"\rho_j"]
    &&Dj'
      \ar[d,"\rho_{j'}"]
  \\ERj
      \ar[rr,"ERh"]
      \ar[dr,swap,near start,"Ef"]
    &&ERj'
      \ar[dl,near start,"Ef'"]
  \\&Ek
  \end{tikzcd}
\]
commutes;
however, the fact that the diagram
\[
  % https://q.uiver.app/#q=WzAsNSxbMCwwLCJEaiJdLFswLDEsIkVSaiJdLFsyLDAsIkVSaiJdLFsyLDEsIkVSaiciXSxbMSwyLCJFayJdLFswLDEsIlxccmhvX2oiLDJdLFswLDIsIlxccmhvX2oiXSxbMiwzLCJFUmgiXSxbMSw0LCJFZiIsMl0sWzMsNCwiRWYnIl1d
  \begin{tikzcd}[column sep=tiny]
    Dj && ERj \\
    ERj && {ERj'} \\
    & Ek
    \arrow["{\rho_j}"', from=1-1, to=2-1]
    \arrow["{\rho_j}", from=1-1, to=1-3]
    \arrow["ERh", from=1-3, to=2-3]
    \arrow["Ef"', from=2-1, to=3-2]
    \arrow["{Ef'}", from=2-3, to=3-2]
  \end{tikzcd}
\]
commutes (which is the case if $h$ defines a morphism in $(R,\rho)/k$) does \emph{not} imply that
\[
  \begin{tikzcd}[column sep=tiny]
    Rj
      \ar[rr,"Rh"]
      \ar[dr,swap,near start,"f"]
    &&Rj'
      \ar[dl,near start,"f'"]
  \\&k
  \end{tikzcd}
\]
commutes. The latter holds only under futher assumptions, such as $E$ being faithful and $\rho_j$ being an epimorphism.

\begin{definition}
\label{definition:relatively-initial}
  Let $(R,\rho)\colon(\cat{J},D)\to(\cat{K},E)$ be a morphism in $\Diag(\cat{X})$.
  We say that the functor $R$ is \emph{initial relative to} the transformation $\rho$ (or that the diagram morphism $(R,\rho)$ is \emph{relatively initial}) if, for all $k\in\cat{K}$, the relative comma category $(R,\rho)/k$ is non-empty and connected.
\end{definition}

Relative initiality \emph{is} a proper generalisation of initiality: a functor $R\colon\cat{J}\to\cat{K}$ is initial if and only if the diagram morphism $(R,\id_R)\colon(\cat{J},R)\to(\cat{K},\id_{\cat{K}})$ is relatively initial, whereas there exist morphisms $(R,\rho)$ that are relatively initial but where the functor $R$ is not initial, such as in \cref{example:relatively-initial-but-not-initial}.
A relatively initial morphism is still sufficient to give a weak equivalence, as the following lemma shows.

\begin{lemma}
\label{lemma:relatively-initial-satisfies-pushforward}
  Let $(R,\rho)$ be a relatively initial morphism in $\Diag(\cat{X})$.
  Then $(R,\rho)$ is a weak equivalence.
\end{lemma}

\begin{proof}
  This is \cite[Theorem~9.12]{patterson2022}.
\end{proof}

Relatively initial morphisms satisfy some standard useful properties: they are closed under composition, and every isomorphism is relatively initial.
Furthermore, relatively initial morphisms preserve limits, in much the same sense that initial functors do.
For more details, see \cite[Proposition~9.13 and Lemma~9.14]{patterson2022}.

\begin{example}
\label{example:relatively-initial-but-not-initial}
  Let $\cat{J}$ be the walking arrow, or the category with two objects and one non-identity morphism, and let $\cat{K}$ be the category with two objects and two parallel non-identity morphisms:
  \[
    \cat{J} = \left\{0 \overset{\alpha}{\rightarrow} 1\right\}
    \qquad\text{and}\qquad
    \cat{K} = \left\{0 \underset{\beta}{\overset{\alpha}{\rightrightarrows}} 1\right\}.
  \]
  Let $f\colon x\to y$ be a morphism in some category $\cat{X}$, and define diagrams $D\colon\cat{J}\to\cat{X}$ and $E\colon\cat{K}\to\cat{X}$ by
  \[
    D(\cat{J}) = \left\{x \xrightarrow{f} y\right\}
    \qquad\text{and}\qquad
    E(\cat{K}) = \left\{x \underset{f}{\overset{f}{\rightrightarrows}} y\right\}.
  \]

  If we define $R\colon\cat{J}\to\cat{K}$ to be the identity-on-objects functor implied by the naming of the morphisms (i.e. $R$ sends the unique morphism $\alpha\colon0\to1$ in $\cat{J}$ to the morphism $\alpha\colon0\to1$ in $\cat{J}$), then this defines a strict morphism of diagrams $(R,\id)\colon(\cat{J},D)\to(\cat{K},E)$ in $\Diag(\cat{X})$.
  Then the following three things are true:
  \begin{enumerate}[noitemsep]
    \item $R$ \emph{is not} initial; but
    \item $(R,\id)$ \emph{is} a weak equivalence; and
    \item $R$ \emph{is} initial with respect to $\id$.
  \end{enumerate}

  The second point follows from the third point combined with \cref{lemma:relatively-initial-satisfies-pushforward}, so it remains only to justify the two concerning initiality.

  That $R$ is not initial follows by considering the comma category $R/1$, which is non-empty but not connected: there is no morphism $\sigma\colon0\to0$ in $\cat{K}$ such that $\alpha\sigma=\beta$ or $\alpha=\beta\sigma$, so we can never hope to find a zigzag between $\alpha\colon R0\to1$ and $\beta\colon R0\to1$.

  However, $(R,\id)$ is \emph{relatively} initial, since the identity morphism $\id\colon0\to0$ in $\cat{J}$ gives a morphism between $\alpha\colon R0\to1$ and $\beta\colon R0\to1$, since the diagram
  \[
    \begin{tikzcd}[column sep=tiny]
      D0
        \ar[rr,"D\id"]
        \ar[d,swap,"\id"]
      &&D0
        \ar[d,"\rho_0"]
    \\ER0
        \ar[dr,swap,"E\alpha"]
      &&ER0
        \ar[dl,"E\beta"]
    \\&E1
    \end{tikzcd}
    \quad=\quad
    \begin{tikzcd}[column sep=tiny]
      X
        \ar[rr,"\id"]
        \ar[d,swap,"\id"]
      &&X
        \ar[d,"\id"]
    \\X
        \ar[dr,swap,"f"]
      &&X
        \ar[dl,"f"]
    \\&Y
    \end{tikzcd}
  \]
  commutes.
  This example illustrates how relatively initiality is ``initiality evaluated in $\cat{X}$.''
\end{example}

Importantly, although relative initiality is strictly more general than initiality and is still a sufficient condition for a morphism of diagrams to be a weak equivalence, it is also still \emph{not necessary}.

\begin{lemma}
\label{lemma:relatively-initial-but-not-pushforward}
  Let $(R,\rho)$ be a morphism in $\Diag(*)$, where $*$ denotes the terminal category.
  Then $(R,\rho)$ is a weak equivalence if and only if $R$ induces a bijection on connected components, and is relatively initial if and only if it \emph{further} satisfies the condition that every $R/k$ is non-empty.
\end{lemma}

\begin{proof}
  Since $*$ is the terminal category, there is an isomorphism $\Diag(*)\cong\Cat$. So, a diagram $D\colon\cat{J}\to*$ is just the information of the indexing category $\cat{J}$, and we can view any functor $R\colon\cat{J}\to\cat{K}$ as a morphism $(R,\id)$ in $\Diag(*)$.

  Now, a discrete opfibration over $*$ is precisely a discrete category, i.e., a set.
  This means that, given a discrete opfibration $\pi\colon\cat{E}\to*$, a lift of $D$ along $\pi$ is simply an arbitrary function of sets $\pi_0\cat{J}\to\cat{E}$, where $\pi_0\cat{J}$ denotes the set of connected components of $\cat{J}$.
  So asking for a morphism $(R,\id)\colon\cat{J}\to\cat{K}$ to be a weak equivalence is thus precisely asking that, for any set $\cat{E}$ and any function $\overline{D}\colon\pi_0\cat{J}\to\cat{E}$, there exists a unique function $\overline{E}\colon\pi_0\cat{K}\to\cat{E}$ such that $\overline{D}=\overline{E}\circ\pi_0R$.
  In particular, this must be true for the set $\cat{E}=\pi_0\cat{J}$ and the function $\overline{D}=\id_{\pi_0\cat{J}}$.
  But this then says that $\pi_0R: \pi_0\cat{J}\to\pi_0\cat{K}$ is a split monomorphism with unique splitting, and thus $\pi_0R$ is a bijection.
  Conversely, if $\pi_0R$ is a bijection, then $(R,\id)$ is a weak equivalence: given any function $\overline{D}\colon\pi_0\cat{J}\to\cat{E}$, we obtain a lift $\overline{E}\colon\pi_0\cat{K}\to\cat{E}$ simply by $\overline{D}$ with the inverse of $\pi_0R$.

  Now consider what it means for a morphism $(R,\id)\colon\cat{J}\to\cat{K}$ to be relatively initial: the relative comma categories $(R,\id)/k$ need to be non-empty and connected for all $k\in\cat{K}$.
  But the morphisms in these relative comma categories are just the morphisms in $\cat{J}$, since the diagram in \cref{definition:relative-slice-category} consists of objects and morphisms in $\cat{X}$, which is here the trivial category $*$ and so this diagram commutes for any choice of morphism in $\cat{J}$.
  In other words, $(j,f)$ and $(j',f')$ are in the same connected component of $(R,\id)/k$ if and only if $j$ and $j'$ are in the same connected component of $\cat{J}$.
  So if $\pi_0$ is a bijection \emph{and} every $(R,\id)/k$ is non-empty, then $R$ is relatively initial.
  Conversely, if $R$ is relatively initial then every $(R,\id)/k$ is non-empty by definition, and thus $\pi_0R$ is surjective;
  but $\pi_0R$ is also injective, since if $j$ and $j'$ are in the same connected component of $\cat{J}$ then $Rj$ and $Rj'$ are in the same connected component of $\cat{K}$.
\end{proof}

\section{Localisations of categories of diagrams}
\label{section:localisations-of-categories-of-diagrams}

Motivated by the idea that diagrams present systems of equations, we have so far
proposed a notion of equivalence between diagrams, which we call a weak
equivalence in $\DiagOp(\cat{X})$ (\cref{definition:equivalence-of-systems}); we showed
that this is essentially equivalent to being a weak equivalence in $\Diag(\cat{X})$, at least for pseudo morphisms
(\cref{corollary:pushforward-iff-equivalence}); and we saw that initiality and,
more generally, relative initiality are sufficient
(\cref{lemma:relatively-initial-satisfies-pushforward}) but not necessary
(\cref{example:relatively-initial-but-not-initial} and
\cref{lemma:relatively-initial-but-not-pushforward}) conditions to have a weak
equivalence.

Now, if we want to take seriously the consideration that two diagrams
in the category $\cat{X}$ should be ``the same'' whenever they are weakly equivalent, then we
need to formally invert the weak equivalences, i.e. freely turn them into isomorphisms. The tool to achieve this is known as
\emph{localisation} \cite{gabriel1967}.

\begin{definition}
\label{definition:localisation}
  Let $\cat{W}$ be a class of morphisms in a category $\cat{C}$.
  The \emph{localisation of $\cat{C}$ along $\cat{W}$} is, if it exists, the data of a category $\cat{C}[\cat{W}^{-1}]$ along with a functor $Q\colon\cat{C}\to\cat{C}[\cat{W}^{-1}]$ such that
  \begin{enumerate}[i.]
    \item for all morphisms $f\in\cat{W}$, the morphism $Qf$ in $\cat{C}[\cat{W}^{-1}]$ is an isomorphism;
    \item any functor $F\colon\cat{C}\to\cat{D}$ such that $Ff$ is an isomorphism for all $f\in W$ factors uniquely through $Q\colon\cat{C}\to\cat{C}[\cat{W}^{-1}]$.
  \end{enumerate}
  The universal property of the localisation makes $\cat{C}[\cat{W}^{-1}]$ the initial category in which all morphisms in $\cat{W}$ become invertible.
\end{definition}

The localisation is in general very complicated. One condition on the class of weak equivalences that
begins to make the localisation more manageable without being too stringent is as follows.

\begin{definition}
\label{definition:2-out-of-3}
  Let $\cat{C}$ be a category and let $\cat{W}$ be a class of morphisms in $\cat{C}$.
  We say that $\cat{W}$ satisfies \emph{2-out-of-3} if, for any pair of composable morphisms $f,g$ in $\cat{C}$, if any two of $f$, $g$, and $gf$ are in $\cat{W}$ then so too is the third.
\end{definition}

\begin{example}
\label{example:2-out-of-6}
  The prototypical class of morphisms that satisfies 2-out-of-3 is the class of isomorphisms in any category.

  Since weak equivalences in $\DiagOp$ are defined in terms of families of bijections, they also satisfy 2-out-of-3.
  The same thus holds, by \cref{corollary:pushforward-iff-equivalence}, for weak equivalences in $\strongDiag$.
\end{example}

\begin{remark}
  \label{remark:initial-does-not-satisfy-2-out-of-3}
    \cref{example:relatively-initial-but-not-initial} shows that initial functors do \emph{not} satisfy 2-out-of-3, in contrast to weak equivalences.

    Indeed, consider the unique identity-on-objects functor $S\colon\cat{K}\to\cat{J}$, which sends both $\alpha$ and $\beta$ to $\alpha$.
    Then one can show that $S$ is initial, and $SR=\id_\cat{J}$ is initial since identity functors always are, but we already know that $R$ is not initial. In particular, this shows that the localisation of $(\Cat/\cat{X})$ along the class of initial functors contains
    isomorphisms that do not arise directly from initial functors.
    This is an important
    point to keep in mind.
  \end{remark}

We now return to the specific case of categories of diagrams.
As a preliminary remark, recall that strict diagram categories are isomorphic to slices of $\Cat$:
\[
  \strictDiag(\cat{X})
  \cong \Cat/\cat{X}
  \cong \strictDiagOp(\cat{X})^\op.
\]
Thus, we can talk about morphisms in $\Cat/\cat{X}$ as weak equivalences in $\Diag(\cat{X})$ or $\DiagOp(\cat{X})$ by regarding them as morphisms in the relevant diagram category.

\begin{definition}
  Write $\iclass$ for the class of morphisms in $\Cat/\cat{X}$ given by initial functors, $\pclass$ for the class of weak equivalences in $\Diag(\cat{X})$, and $\eosclass$ for the class of weak equivalences in $\DiagOp(\cat{X})$.

  We also write $\pclass^\strict\subset\pclass$ and $\eosclass^\strict\subset\eosclass$ for the subclasses of strict morphisms.
  By the remark above, we can regard $\pclass^\strict$ and $\eosclass^\strict$ as classes of weak equivalences in $\Cat/\cat{X}$ and $(\Cat/\cat{X})^\op$, respectively.
\end{definition}

Our original motivation suggests that we are interested in the category $\DiagOp(\cat{X})[\eosclass^{-1}]$. We know that initiality is a sufficient condition for being a weak equivalences, so this category admits a canonical map from $(\Cat/\cat{X})[\iclass^{-1}]$.
The main theorem of this paper says that, in fact, this
map is an isomorphism.

\begin{theorem}
\label{theorem:localisation-of-slice-cat-and-of-diag-op}
  The inclusion of $\Cat/\cat{X}^{\op}$ into $\DiagOp(\cat{X})$ induces an isomorphism
  \[
    (\Cat/\cat{X})^{\op}[\iclass^{-1}]
    \xrightarrow{\cong} \DiagOp(\cat{X})[\eosclass^{-1}]
  \]
  from the localisation of $\Cat/\cat{X}^\op$ at the class $\iclass$ of initial functors to the localisation of $\DiagOp(\cat{X})$ at the class $\eosclass$ of weak equivalences.
\end{theorem}

\begin{remark}
	Some remarks on the nature of the category on the left-hand side of the isomorphism in \cref{theorem:localisation-of-slice-cat-and-of-diag-op} are in order. First, since $\Cat/\cat{X}$ admits
	a lifted comprehensive factorization system, where, as in $\Cat$, the left part consists of initial functors and the right part of discrete opfibrations,
	we can express $\Cat/\cat{X}[\iclass^{-1}]$ more simply as the category of discrete opfibrations over $\cat{X}.$ Indeed, any diagram
	$D\colon\cat J\to \cat{X}$ factors as $\pi\circ i$, where $\pi$ is a discrete opfibration and $i$ is initial, which means that once $i$ is inverted,
	$D$ becomes isomorphic to the discrete opfibration $\pi.$ Furthermore, any morphism in $\Cat/\cat{X}$ between discrete opfibrations is itself
	a discrete opfibration: if $D,D'\colon\cat J,\cat{J'}\to \cat{X}$ are discrete opfibrations and $R:\cat J\to \cat{J'}$ satisfies $D'\circ R=D$, then
	if $R$ were not itself a discrete opfibration then factoring the initial part out of $R$ would give a factorization of $D$ with nontrivial initial
	part, which contradicts the essential uniqueness of the comprehensive factorization.

	The conclusion is that $\Cat/\cat{X}[\iclass^{-1}]$ is equivalent to the category $\cat{Dopf}/\cat{X}$, where $\cat{Dopf}$ is the wide subcategory of
	$\cat{Cat}$ spanned by the discrete opfibrations. It is well-known that, via the Grothendieck construction, $\cat{Dopf}/\cat{X}$ is itself equivalent
	to the category of copresheaves $\cat{X}\to \Set.$ Thus the category of systems of equations in $\cat{X}$ (localised along initial functors) can also be viewed simply as the category
	of copresheaves on $\cat{X}$; in particular, this category forms a topos, which demonstrates the exceptionally rich structure emerging from the
	categorification process relative to the traditional conception of equations.

	That said, we do not prefer to think of $\Cat/\cat{X}[\iclass^{-1}]$ as the category of copresheaves on $\cat{X}$, since the objects of this category
	are difficult to specify in a finitary way in the case where $\cat{X}$ is not itself finitely presented (which is our baseline). Instead, we prefer to
	consider how to express $\Cat/\cat{X}[\iclass^{-1}]$ most simply without changing the objects, since diagrams in $\cat{X}$ are much more amenable to
	finite specification in practical cases. To that end, we observe that the hom-set $\Cat/\cat{X}[\iclass^{-1}](D\colon\cat J\to \cat{X},D'\colon\cat{J'}\to\cat{X})$ between any two diagrams
	is in bijection with the hom-set $\cat{Dopf}/\cat{X}(\pi\colon\bar{\cat J}\to \cat{X},\pi'\colon\overline{\cat{J'}}\to \cat{X})$, where
	$\pi$ and$ \pi'$ are the discrete opfibrations generated by $D$ and $D'$ (respectively) via the comprehensive factorization system. This bijection is naturally induced
	by composition with $i,i'\colon\cat J\to \bar{\cat J},\cat{J'}\to \overline{\cat{J'}}$ respectively, which are the initial parts of the same comprehensive factorisations.
	Thus, each morphism in $\Cat/\cat{X}[\iclass^{-1}](D\colon\cat J\to \cat{X},D'\colon\cat{J'}\to\cat{X})$ may be canonically written as the composite of $i$,
	a discrete opfibration over $\cat{X}$, and $(i')^{-1}$; composing the first two morphisms in $\Cat/\cat{X}$, we obtain a unique expression of such a morphism
	as a zigzag $D\to \pi'\leftarrow D'\colon i'.$ Again, we shall often prefer not to calculate $\pi'$ explicitly, so that the real import of this argument
	is that every morphism in $\Cat/\cat{X}[\iclass^{-1}](D\colon\cat J\to \cat{X},D'\colon\cat{J'}\to\cat{X})$ may be written in \emph{some} manner as a zigzag
	$D\to D''\leftarrow D'$, where the reversed arrow is determined by an initial functor between the domains of $D''$ and $D'.$ It is this expression that
	we expect to find most useful for applications.
\end{remark}

The rest of this section will be dedicated to a proof of \cref{theorem:localisation-of-slice-cat-and-of-diag-op}, which we break down into four steps:
\begin{itemize}
  \item[\emph{Step 1.}] In $\Cat/\cat{X}$, inverting initial functors is the same as inverting weak equivalences:

    $(\Cat/\cat{X})[\iclass^{-1}]\xrightarrow{\cong}(\Cat/\cat{X})[(\pclass^\strict)^{-1}]$
  \item[\emph{Step 2.}] Inverting weak equivalences in $\Cat/\cat{X}$ commutes with taking the opposite category:

    $(\Cat/\cat{X})[(\pclass^\strict)^{-1}]^\op\xleftarrow{\cong}(\Cat/\cat{X})^\op[(\eosclass^\strict)^{-1}]$
  \item[\emph{Step 3.}] Inverting (strict) weak equivalences in $(\Cat/\cat{X})^\op$ is the same as inverting them in all of $\DiagOp(\cat{X})$:

    $(\Cat/\cat{X})^\op[(\eosclass^\strict)^{-1}]\xrightarrow{\cong}\DiagOp(\cat{X})[(\eosclass^\strict)^{-1}]$
  \item[\emph{Step 4.}] In $\DiagOp(\cat{X})$, inverting strict weak equivalences is the same as inverting all weak equivalences:

    $\DiagOp(\cat{X})[(\eosclass^\strict)^{-1}]\xrightarrow{\cong}\DiagOp(\cat{X})[\eosclass^{-1}]$.
\end{itemize}
The composition of these four isomorphisms will give us a proof of \cref{theorem:localisation-of-slice-cat-and-of-diag-op}.

\subsection*{Step 1. Initial functors and weak equivalences in $\Cat/\cat{X}$}

The key technical observation for this step is the following.

\begin{lemma}
\label{lemma:pushforward-between-disc-opfibs-implies-iso}
  Let $P:\cat E\to \cat{X}$ and $Q:\cat F\to \cat{X}$ be discrete opfibrations.
  Suppose we have $R:\cat E\to \cat F$ such that $R\colon(\cat{E},P)\to(\cat{F},Q)$ is a weak equivalence in $\Cat/\cat{X}$. Then $R$ is an isomorphism of categories.
\end{lemma}

\begin{proof}
  Since $R$ is a weak equivalence, it induces a bijection $R_*\colon\Lift{P}{\pi}\to \Lift{Q}{\pi}$ for every discrete opfibration $\pi$ over $\cat{X}.$ In particular we
  may consider $\pi=Q$ and consider the lift $\id_{\cat F}$ of $Q$ through
  itself, shown below:
  \begin{equation*}
    % https://q.uiver.app/#q=WzAsNCxbMCwwLCJcXGNhdCBFIl0sWzEsMCwiXFxjYXQgRiJdLFsxLDEsIlxcY2F0IFgiXSxbMCwxLCJcXGNhdCBGIl0sWzMsMiwiUSIsMl0sWzEsMiwiUSJdLFszLDEsIlxcaWRfe1xcY2F0IEZ9Il0sWzAsMSwiUiJdLFszLDAsIlMiLDAseyJzdHlsZSI6eyJib2R5Ijp7Im5hbWUiOiJkb3R0ZWQifX19XV0=
    \begin{tikzcd}
      {\cat E} & {\cat F} \\
      {\cat F} & {\cat{X}}
      \arrow["Q"', from=2-1, to=2-2]
      \arrow["Q", from=1-2, to=2-2]
      \arrow["{\id_{\cat F}}", from=2-1, to=1-2]
      \arrow["R", from=1-1, to=1-2]
      \arrow["S", dotted, from=2-1, to=1-1]
    \end{tikzcd}
  \end{equation*}
Then, since $R_*$ is a bijection, we can find a unique $S$ such that $RS=\id_{\cat F}$
(and, in this case redundantly, $PR=Q.$) Thus $R$ is a split epimorphism and
$S$ is a split monomorphism. Now, since $R$ and $\id_{\cat F}$ are weak equivalences,
so is $S$ by two-out-of-three, so that we can repeat the same argument with $S$ in
place of $R$ to show that $S$ is also a split epimorphism. Therefore, $S$ is an
isomorphism, and finally so is $R.$
\end{proof}

We can now apply this lemma to prove the following.

\begin{lemma}
\label{lemma:localisation-at-initial-localises-pushforwards}
  Let $L_\iclass\colon\Cat/\cat{X}\to(\Cat/\cat{X})[\iclass^{-1}]$ be the localisation functor that inverts the class of initial functors.
  Then $L_\iclass$ sends weak equivalences to isomorphisms.
\end{lemma}

\begin{proof}
  Let $R\colon(\cat{J},D)\to(\cat{K},E)$ be a weak equivalence in $\Cat/\cat{X}$.

  Using the comprehensive factorisation system (\cref{lemma:comprehensive-factorisation-system}), we can factor $E$ as $PI$, with $I\colon\cat{K}\to\cat{L}$ initial and $P\colon\cat{L}\to\cat{X}$ a discrete opfibration.
  This gives us the commutative diagram
  \[
    \begin{tikzcd}
      \cat{J}
        \ar[r,"R"]
        \ar[dr,swap,"D"]
      &\cat{K}
        \ar[r,"I"]
        \ar[d,"E"]
      &\cat{L}
        \ar[dl,"P"]
    \\&\cat{X}
    \end{tikzcd}
  \]
  Since $I$ is initial, it is a weak equivalence, by \cref{lemma:initial-satisfies-pushforward}, and $R$ is a weak equivalence by assumption;
  thus the composite $IR$ is also a weak equivalence.
  It thus suffices to show that the composite $IR$ is initial, since then $(IR)^{-1}I$ will be an inverse for $R$ in $(\Cat/\cat{X})[\iclass^{-1}]$.

  If we again apply the comprehensive factorisation system to factor $D\colon\cat{J}\to\cat{X}$ as $QJ$, with $J\colon \cat J \to \cat M$ initial and $Q\colon\cat M \to \cat{X}$ a discrete opfibration, then we obtain the commutative square
  \[
    \begin{tikzcd}
      \cat{J}
        \ar[r,"IR"]
        \ar[d,swap,"J"]
      &\cat{L}
        \ar[d,"P"]
    \\\cat{M}
        \ar[r,swap,"Q"]
      &\cat{X}
    \end{tikzcd}.
  \]
  The unique lifting property of orthogonal factorisation systems tells us that the square above has a unique diagonal filler $F\colon\cat{M}\to\cat{L}$.
  In particular, we find that $IR=FJ$ as morphisms $D\to P$ in $\Cat/\cat{X}$:
  \[
    \begin{tikzcd}[column sep=tiny]
      \cat{J}
        \ar[d,swap,"J"]
        \ar[rr,"R"]
        \ar[ddr,bend right=80,swap,"D"]
      &&\cat{K}
        \ar[d,"I"]
        \ar[ddl,bend left=80,"E"]
    \\\cat{M}
        \ar[rr,"F"]
        \ar[dr,swap,"Q"]
      &&\cat{L}
        \ar[dl,"P"]
    \\&\cat{X}
    \end{tikzcd}.
  \]
  Finally, since $IR$ is a weak equivalence, and since $J$ is initial and is thus a weak equivalence, applying 2-out-of-3 again tells us $F$ is also a weak equivalence.
  We can now apply \cref{lemma:pushforward-between-disc-opfibs-implies-iso}, which tells us that $F\colon(\cat{M},Q)\to(\cat{L},P)$, as a weak equivalence between two discrete opfibrations, is an isomorphism.
  But then $IR=FJ$ is the composite of an initial functor $J$ with an isomorphism $F$, and so is also initial, as was to be shown.
\end{proof}

\begin{corollary}
\label{corollary:localisation-step-1}
  The canonical functor
  \[
    (\Cat/\cat{X})[\iclass^{-1}]
    \to (\Cat/\cat{X})[(\pclass^\strict)^{-1}]
  \]
  is an isomorphism.
\end{corollary}

\begin{proof}
  We know from \cref{lemma:initial-satisfies-pushforward} that $\iclass\subseteq\pclass^\strict$, so inverting $\pclass^\strict$ will in particular invert $\iclass$.
  Conversely, \cref{lemma:localisation-at-initial-localises-pushforwards} tells us that inverting $\iclass$ also inverts $\pclass^\strict$ in the process.
  All together then, inverting $\iclass$ is equivalent to inverting $\pclass^\strict$.
\end{proof}

\subsection*{Step 2. Weak equivalences in $\Cat/\cat{X}$ and $(\Cat/\cat{X})^\op$}

In \cref{corollary:pushforward-iff-equivalence}, we showed that taking mates induces a correspondence between weak equivalences in $\strongDiag(\cat{X})$ and weak equivalences in $\strongDiagOp(\cat{X})$.
In fact, the mate actually gives an isomorphism of categories
\[
  \strongDiag(\cat{X})
  \cong \strongDiagOp(\cat{X})^\op.
\]
\Cref{corollary:pushforward-iff-equivalence} then tells us that this induces an isomorphism of wide subcategories between pseudo morphisms that are weak equivalences in $\strongDiag(\cat{X})$ and pseudo morphisms that are weak equivalences in $\strongDiagOp(\cat{X})$.

It follows immediately from the definition of localisation of categories that an isomorphism $(\cat{C},\cat{A})\cong(\cat{D},\cat{B})$ of pairs of a category and a wide subcategory induces an isomorphism $\cat{C}[\cat{A}^{-1}]\cong\cat{D}[\cat{B}^{-1}]$ of localisations.
Restricting from pseudo morphisms to strict morphisms, this means that we have proved the following.

\begin{corollary}
\label{corollary:localisation-step-2}
  The canonical functor
  \[
    (\Cat/\cat{X})^\op[(\eosclass^\strict)^{-1}]
    \to (\Cat/\cat{X})[(\pclass^\strict)^{-1}]^\op
  \]
  is an isomorphism.
\end{corollary}

\subsection*{Step 3. Weak equivalences in $\Cat/\cat{X}$, strict weak equivalences in $\DiagOp(\cat{X})$}

The key technical observations for this step are the following lemmas, where $D\colon\cat{J}\to\cat{X}$ is an object of $\Cat/\cat{X}$ throughout.

Given a cospan $F:\cat A\to \cat C\leftarrow \cat B:G$, we shall use
the below notations for the comma category and the components of its canonical
cone:
% https://q.uiver.app/#q=WzAsNCxbMCwxLCJcXGNhdCBBIl0sWzEsMSwiXFxjYXQgQyJdLFsxLDAsIlxcY2F0IEIiXSxbMCwwLCJGL0ciXSxbMCwxLCJGIiwyXSxbMiwxLCJHIl0sWzMsMCwiXFxwaV97Ri99IiwyXSxbMywyLCJcXHBpX3svR30iXSxbNiw1LCJcXGFscGhhX3tGL0d9IiwwLHsic2hvcnRlbiI6eyJzb3VyY2UiOjMwLCJ0YXJnZXQiOjMwfX1dXQ==
\[\begin{tikzcd}
	{F/G} & {\cat B} \\
	{\cat A} & {\cat C}
	\arrow["F"', from=2-1, to=2-2]
	\arrow[""{name=0, anchor=center, inner sep=0}, "G", from=1-2, to=2-2]
	\arrow[""{name=1, anchor=center, inner sep=0}, "{\pi_{F/}}"', from=1-1, to=2-1]
	\arrow["{\pi_{/G}}", from=1-1, to=1-2]
	\arrow["{\omega_{F/G}}", shorten <=11pt, shorten >=11pt, Rightarrow, from=1, to=0]
\end{tikzcd}\]
In the case that $F=\id_{\cat C}$, we'll write $\cat C/G$ and $\pi_{\cat C/}$,
and similarly if $G=\id_{\cat C}.$

For clarity and compatibility with the 2-categorical section later on, we
work mainly with the universal property of the comma category: namely,
that for any cospan $(F,G)$ as above and for any $\cat D$, the hom-category $\Cat(\cat D,F/G)$
is canonically isomorphic to the category whose objects are triples
$(A\colon\cat D\to \cat A,B\colon\cat D\to \cat B,\gamma\colon FA\to GB)$ and whose
morphisms $(A,B,\gamma)\to (A',B',\gamma')$ are pairs of
natural transformations $(\alpha\colon A\to A',\beta\colon B\to B')$ such that
the square below commutes:% https://q.uiver.app/#q=WzAsNCxbMCwwLCJGQSJdLFsxLDAsIkdCIl0sWzAsMSwiRkEnIl0sWzEsMSwiR0InIl0sWzAsMiwiRlxcYWxwaGEiLDJdLFsxLDMsIkdcXGJldGEiXSxbMCwxLCJcXGdhbW1hIl0sWzIsMywiXFxnYW1tYSciLDJdXQ==
\[\begin{tikzcd}
	FA & GB \\
	{FA'} & {GB'}
	\arrow["F\alpha"', from=1-1, to=2-1]
	\arrow["G\beta", from=1-2, to=2-2]
	\arrow["\gamma", from=1-1, to=1-2]
	\arrow["{\gamma'}"', from=2-1, to=2-2]
\end{tikzcd}\]
Specifically, this isomorphism is induced by sending a functor $T\colon\cat D\to F/G$ to the triple $(\pi_{F/}T,\pi_{/G}T,\omega_{F/G}T)$
and a natural transformation $\tau:T\To T':\cat D\to F/G$ to the pair $(\pi_{F/}\tau,\pi_{/G}\tau).$

\begin{lemma}
\label{lemma:slice-into-diag-has-adjoint}
  The canonical inclusion of $\Cat/\cat{X}$ into $\Diag(\cat{X})$ has a right adjoint given on objects by sending $D\colon\cat{J}\to\cat{X}$ to the canonical projection functor $\pi_{\cat{X}/}\colon\cat{X}/D\to\cat{X}.$

  Dually, the canonical inclusion of $(\Cat/\cat{X})^\op$ into $\DiagOp(\cat{X})$ has a left adjoint given on objects by sending $D\colon\cat{J}\to\cat{X}$ to the canonical projection functor $\pi_{/\cat{X}}\colon D/\cat{X}\to\cat{X}.$
\end{lemma}

\begin{proof}
  We want to construct a natural isomorphism
  \[
    \Diag(\cat{X})\left( \cat{J}\xrightarrow{D}\cat{X}, \cat{K}\xrightarrow{E}\cat{X} \right)
    \cong
    \Cat/\cat{X}\left( \cat{J}\xrightarrow{D}\cat{X}, (\cat{X}/E)\xrightarrow{\pi_{\cat{X}/}}\cat{X} \right).
  \]

  By the universal property described just above the statement of the lemma, to
  give a functor $S\colon\cat J\to (\cat{X}/E)$ it is equivalent to specify the functors
  $\pi_{\cat{X}/} S\colon\cat J\to \cat{X}$ and $\pi_{/E} S\colon\cat J\to \cat K$,
  and the natural transformation $\omega_{\cat{X}/E} S\colon \pi_{\cat{X}/} S\To E\pi_{/E} S.$ Such an $S$
  gives a morphism in $\Cat/\cat{X}$ if and only if $\pi_{\cat{X}/} S=D$,
  so that a morphism in $\Cat/\cat{X}(D,\pi_{\cat{X}/})$ is precisely
  given by a functor $R=\pi_{/E}S\colon\cat J\to \cat K$ and a natural transformation
  $D\To ER$, which is the same as the data of a morphism on the left-hand side.

\end{proof}

\begin{lemma}
\label{lemma:projections-from-comma-are-split-adjoints}
  Consider the functor $\iota_{/D}:\cat J\to \cat{X}/D$ determined by
  \[
    \pi_{\cat{X}/} \iota_{/D}= D,\quad \pi_{/D} \iota_{/D} =\id_{\cat J},\quad  \omega_{\cat{X}/D}\iota_{/D} = \id_D.
  \]
  Then $\iota_{/D}$ is right adjoint to $\pi_{/D}.$

  Dually, $\pi_{D/}\colon D/\cat{X}\to\cat{J}$ has a left adjoint $\iota_{D/}$ that is a monomorphism that splits $\pi_{D/}$ and defines a strict morphism of diagrams $(\cat{J},D)\to(D/\cat{X},\pi_{/\cat{X}})$.
\end{lemma}

Before we begin the proof, note that by definition, $\iota_{/D}$ is a split monomorphism splitting $\pi_{/D}$.
Furthermore, $\iota_{/D}$ defines a strict morphism $(\cat{J},D)\to(\cat{X}/D,\pi_{\cat{X}/})$ in $\Cat/\cat{X}$.
\begin{proof}
  Define the counit $\pi_{/D} \iota_{/D}\to \id_{\cat J}$ to be the identity
  of $\id_{\cat J}$, since $\pi_{/D}$ splits $\iota_{/D}.$ Next, the unit must
  be a map $\eta:\id_{\cat{X}/D}\to \iota_{/D} \pi_{/D}.$ First define
  $\pi_{/D}\eta$ as $\id_{\pi_{/D}}.$ For $\pi_{\cat{X}/}\eta$,
  we need a natural transformation $\pi_{\cat{X}/}\to \pi_{\cat{X}/}\iota_{/D}\pi_{/D}=D\pi_{/D}$; we are given such a natural transformation in
  $\omega_{\cat{X}/D}.$ This determines a morphism $\eta$ as desired
  since $\omega_{\cat{X}/D}\id_{\cat{X}/D}=\omega_{\cat{X}/D}$ while
  $\omega_{\cat{X}/D}\iota_{/D}\circ \pi_{/D}=\id_{\pi_{/D}}.$

  The triangle identities will now follow once we check that $\pi_{/D}\eta$
  and $\eta \iota_{/D}$ are isomorphisms, see \cite[B.4.2]{riehl-verity}.
  But we defined $\pi_{/D}\eta$ as $\id_{\pi/D}.$ Meanwhile,
  $\eta \iota_{/D}:\iota_{/D}\To \iota_{/D}\circ \pi_{/D}\circ \iota_{/D}:\cat J\to \cat{X}/D$ will be an isomorphism as soon as its whiskerings with
  $\pi_{/D}$ and $\pi_{\cat{X}/}$ are so. Again, the whiskering with $\pi_{/D}$ is
  invertible by assumption; finally $\pi_{\cat{X}/} \eta \iota_{/D}=\omega_{\cat{X}/D}\iota_{/D}=\id_D$, and the claim is proven.
\end{proof}

Of these two adjoint pairs, $\pi_{/\cat{J}}\dashv\iota_{/\cat{J}}$ and $\iota_{\cat{J}/}\dashv\pi_{\cat{J}/}$, we are interested in the latter.
This is because, although $\iota_{/\cat{J}}$ and $\iota_{\cat{J}/}$ both define strict morphisms of diagrams, it is only the latter which is a \emph{left} adjoint and thus initial (\cref{lemma:left-adjoints-are-initial}).
This allows us to prove the following.

\begin{corollary}
\label{corollary:projection-from-comma-is-equivalence}
  Every canonical projection functor $\pi_{D/}\colon D/\cat{X}\to\cat{J}$ is a weak equivalence, as is its left adjoint $\iota_{D/}\colon\cat{J}\to D/\cat{X}$.
\end{corollary}

\begin{proof}
  Firstly, the fact that $\iota_{D/}$ is a weak equivalence follows from the fact that it is a
  left adjoint, thus initial (\cref{lemma:left-adjoints-are-initial}),
  thus a weak equivalence (\cref{lemma:initial-satisfies-pushforward}.
  Then, to show that $\pi_{D/}$ is a weak equivalence, note that $\pi_{D/}\iota_{D/}=\id_\cat{J}$.
  But $\id_\cat{J}$ is a weak equivalence, so 2-out-of-3 tells us that $\pi_{D/}$ is also a weak equivalence.
\end{proof}

The above lemmas allow us to prove the main result of this step using the following approach.

Consider the commutative diagram
\begin{equation}
\label{equation:diagram-before-P}
  \begin{tikzcd}[row sep=huge]
    (\Cat/\cat{X})^\op
      \ar[r,"i"]
      \ar[d,swap,"L"]
    & \DiagOp(\cat{X})
      \ar[d,"M"]
  \\(\Cat/\cat{X})^\op[(\eosclass^\strict)^{-1}]
      \ar[r,swap,"i'"]
    & \DiagOp(\cat{X})[(\eosclass^\strict)^{-1}]
  \end{tikzcd}
\tag{A}
\end{equation}
where $L$ and $M$ are the canonical localisation functors, $i$ is the inclusion, and $i'$ is the inclusion induced by $i$ and functoriality of localisation.
We want to show that $i'$ is actually an isomorphism of categories, and we will do so by constructing a functor $P\colon\DiagOp(\cat{X})\to(\Cat/\cat{X})^\op[\cat{W}_\strict^{-1}]$ such that the diagram
\begin{equation}
\label{equation:diagram-for-P}
  \begin{tikzcd}[row sep=huge]
    (\Cat/\cat{X})^\op
      \ar[r,"i"]
      \ar[d,swap,"L"]
    & \DiagOp(\cat{X})
      \ar[d,"M"]
      \ar[dl,"P"description]
  \\(\Cat/\cat{X})^\op[(\eosclass^\strict)^{-1}]
      \ar[r,swap,"i'"]
    & \DiagOp(\cat{X})[(\eosclass^\strict)^{-1}]
  \end{tikzcd}
\tag{B}
\end{equation}
commutes.
We will then show how this implies that $P$ induces a functor $P'$ that is inverse to $i'$, i.e. such that
\begin{equation}
\label{equation:diagram-for-P'}
  \begin{tikzcd}[row sep=huge]
    (\Cat/\cat{X})^\op
      \ar[r,"i"]
      \ar[d,swap,"L"]
    & \DiagOp(\cat{X})
      \ar[d,"M"]
      \ar[dl,"P"description]
  \\(\Cat/\cat{X})^\op[(\eosclass^\strict)^{-1}]
      \ar[r,shift right,swap,"i'"]
    & \DiagOp(\cat{X})[(\eosclass^\strict)^{-1}]
      \ar[l,shift right,swap,"P'"]
  \end{tikzcd}
\tag{C}
\end{equation}
also commutes.

\begin{lemma}
\label{lemma:diagram-for-P-commutes}
  There exists a functor $P\colon\DiagOp(\cat{X})\to(\Cat/\cat{X})^\op[(\eosclass^\strict)^{-1}]$ such that Diagram~\labelcref{equation:diagram-for-P} commutes.
\end{lemma}

\begin{proof}
  We proceed in three steps: first, we construct $P$ as a morphism of graphs (i.e. without checking for functoriality); second, we show that $P$ is indeed functorial; finally, we show that $P$ does indeed make Diagram~\labelcref{equation:diagram-for-P} commute. The main idea of the proof is fairly simple: define $P$ using the canonical factorisation given by the unit of the adjunction $\iota_{\cat{J}/}\dashv\pi_{\cat{J}/}$ from \cref{lemma:slice-into-diag-has-adjoint}.

  \emph{Step 1: $P$ as a graph morphism}

  On \emph{objects}, we define $P\colon\DiagOp(\cat{X})\to(\Cat/\cat{X})^\op[(\eosclass^\strict)^{-1}]$ to be the identity, sending $D\colon\cat{J}\to\cat{X}$ to itself.
  Now, given a \emph{morphism} $(R,\rho)\colon(\cat{J},D)\to(\cat{K},E)$ in $\DiagOp(\cat{X})$, we want to define the associated morphism $P(R,\rho)$
  in $\Cat/\cat{X}.$ As promised, we start by factoring $(R,\rho)$ through the unit of the adjunction $\iota_{\cat{J}/}\dashv\pi_{\cat{J}/}$ from \cref{lemma:slice-into-diag-has-adjoint}, via the equality
  \[
    \begin{tikzcd}
    % https://q.uiver.app/#q=WzAsMyxbMCwwLCJcXGNhdHtKfSJdLFsyLDAsIlxcY2F0e0t9Il0sWzEsMSwiXFxjYXR7WH0iXSxbMSwwLCJSIiwyXSxbMCwyLCJEIiwyXSxbMSwyLCJFIl0sWzQsNSwiXFxyaG8iLDEseyJvZmZzZXQiOi0yLCJzaG9ydGVuIjp7InNvdXJjZSI6MjAsInRhcmdldCI6MjB9fV1d
      {\cat{J}} && {\cat{K}} \\
      & {\cat{X}}
      \arrow["R"', from=1-3, to=1-1]
      \arrow[""{name=0, anchor=center, inner sep=0}, "D"', from=1-1, to=2-2]
      \arrow[""{name=1, anchor=center, inner sep=0}, "E", from=1-3, to=2-2]
      \arrow["\rho"{description}, shift left=2, shorten <=7pt, shorten >=7pt, Rightarrow, from=0, to=1]
    \end{tikzcd}
    \quad=\quad
    \begin{tikzcd}
    % https://q.uiver.app/#q=WzAsNCxbMCwwLCJcXGNhdHtKfSJdLFsyLDAsIlxcY2F0e0t9Il0sWzEsMSwiXFxjYXR7WH0iXSxbMSwwLCJEL1xcY2F0e1h9Il0sWzAsMiwiRCIsMl0sWzEsMiwiRSJdLFszLDAsIlxccGlfe1xcY2F0e0p9L30iLDJdLFsxLDMsIlxcd2lkZWhhdHtcXHJob30iLDJdLFszLDIsIlxccGlfey9cXGNhdHtYfX0iXSxbNCw4LCJcXG9tZWdhIiwwLHsib2Zmc2V0IjotMiwic2hvcnRlbiI6eyJzb3VyY2UiOjIwLCJ0YXJnZXQiOjIwfX1dXQ==
      {\cat{J}} & {D/\cat{X}} & {\cat{K}} \\
      & {\cat{X}}
      \arrow[""{name=0, anchor=center, inner sep=0}, "D"', from=1-1, to=2-2]
      \arrow["E", from=1-3, to=2-2]
      \arrow["{\pi_{D/}}"', from=1-2, to=1-1]
      \arrow["\widehat{\rho}"', from=1-3, to=1-2]
      \arrow[""{name=1, anchor=center, inner sep=0}, "{\pi_{/\cat{X}}}", from=1-2, to=2-2]
      \arrow["\omega", shift left=2, shorten <=4pt, shorten >=4pt, Rightarrow, from=0, to=1]
    \end{tikzcd}
  \]
  Since $\iota_{D/}$ defines a strict morphism $D\leftarrow \pi_{/\cat{X}}$ (\cref{lemma:projections-from-comma-are-split-adjoints})
  that is also a weak equivalence (\cref{corollary:projection-from-comma-is-equivalence}), it is in
  $\eosclass^\strict$ and thus inverted by $M$.
  The fact that $\iota_{D/}$ is also split by $\pi_{D/}$ (\cref{lemma:projections-from-comma-are-split-adjoints}) then implies that $M(\iota_{D/},\id)^{-1}=M(\pi_{D/},\omega)$, since splittings are preserved by all functors and isomorphisms can be split only by their inverse.

  Putting this all together and using also the facts that the image of $i$ consists of the strict morphisms in $\DiagOp(\cat{X})$ and that Diagram~\labelcref{equation:diagram-before-P} commutes, we see that\footnote{Recall that we are working in $\DiagOp(\cat{X})$, so although $R=\pi_{D/}\circ\widehat{\rho}$ as a \emph{functor}, as a \emph{diagram morphism} we write $(R,\rho)=(\widehat{\rho},\id)\circ(\pi_{D/},\omega)$.}
  \[
    \begin{aligned}
      M(R,\rho)
      &= M(\widehat{\rho},\id) M(\pi_{D/},\omega)
    \\&= M(\widehat{\rho},\id) M(\iota_{D/},\id)^{-1}
    \\&= Mi(\widehat{\rho}) Mi(\iota_{D/})^{-1}
    \\&= i'L(\widehat{\rho}) i'L(\iota_{D/})^{-1}.
    \end{aligned}
  \]

  With this in mind, if $i'$ is really to turn out as an isomorphism with $M=i' P$, we are forced to finish the definition of $P$ as:
  \[
    \begin{aligned}
      P\colon \DiagOp(\cat{X})
      &\to (\Cat/\cat{X})^\op[(\eosclass^\strict)^\op]
    \\(\cat{J},D)
      &\mapsto (\cat{J},D)
    \\(R,\rho)
      &\mapsto L(\widehat{\rho}) L(\iota_{D/})^{-1}.
    \end{aligned}
  \]

  \medskip
  \emph{Step 2: Checking that $P$ is a functor}

  When $(R,\rho)=(\id_{\cat J},\id_D)$ is an identity morphism in $\DiagOp(\cat
  X)$, the associated functor $\hat\rho$ is simply $\iota_{\cat J/}$ itself,
  which shows that $P$ preserves identities.

  Now consider a composable pair of morphisms in $\DiagOp(\cat{X})$, as displayed below.
  \[
    \begin{tikzcd}
    % https://q.uiver.app/#q=WzAsNCxbMSwwLCJcXGNhdHtLfSJdLFsyLDAsIlxcY2F0e0x9Il0sWzAsMCwiXFxjYXR7Sn0iXSxbMSwxLCJcXGNhdHtYfSJdLFsyLDMsIkQiLDJdLFswLDIsIlIiLDJdLFsxLDAsIlMiLDJdLFsxLDMsIkYiXSxbMCwzLCJFIiwxXSxbNCw4LCJcXHJobyIsMCx7Im9mZnNldCI6LTIsInNob3J0ZW4iOnsic291cmNlIjoyMCwidGFyZ2V0IjoyMH19XSxbOCw3LCJcXHNpZ21hIiwwLHsib2Zmc2V0IjotMiwic2hvcnRlbiI6eyJzb3VyY2UiOjIwLCJ0YXJnZXQiOjIwfX1dXQ==
      \cat{J} & \cat{K} & \cat{L} \\
      & \cat{X}
      \arrow[""{name=0, anchor=center, inner sep=0}, "D"', from=1-1, to=2-2]
      \arrow["R"', from=1-2, to=1-1]
      \arrow["S"', from=1-3, to=1-2]
      \arrow[""{name=1, anchor=center, inner sep=0}, "F", from=1-3, to=2-2]
      \arrow[""{name=2, anchor=center, inner sep=0}, "E"{description}, from=1-2, to=2-2]
      \arrow["\rho", shift left=2, shorten <=3pt, shorten >=3pt, Rightarrow, from=0, to=2]
      \arrow["\sigma", shift left=2, shorten <=3pt, shorten >=3pt, Rightarrow, from=2, to=1]
    \end{tikzcd}
  \]
  By definition, we must compare the two right-hand sides below:
  \[
    \begin{aligned}
      P(S,\sigma) P(R,\rho)
      &&&= L(\widehat{\sigma}) L(\iota_{E/})^{-1}  L(\widehat{\rho})  L(\iota_{D/})^{-1}
    \\P((S,\sigma)(R,\rho))
      &= P(RS,\sigma\rho)
    &&= L(\widehat{\sigma\rho})\circ L(\iota_{D/})^{-1}.
    \end{aligned}
  \]
  Cancelling the isomorphism $L(\iota_{D/})^{-1}$, it suffices to show that
  \[
    L(\widehat{\sigma}) L(\iota_{E/})^{-1}  L(\widehat{\rho})
    = L(\widehat{\sigma\rho}).
  \]
  Working now in $\Cat/\cat{X}$, rather than its opposite, we have the four functors
  \[
    \begin{tikzcd}[sep=small]
      & E/\cat{X}
    \\\cat{K}
        \ar[ur,"\iota_{E/}"]
        \ar[dr,swap,"\widehat{\rho}"]
      && \cat{L}
        \ar[ul,swap,"\widehat{\sigma}"]
        \ar[dl,"\widehat{\sigma\rho}"]
    \\& D/\cat{X}
    \end{tikzcd}
  \]
  and we want to show that the two paths from $\cat{L}$ to $D/\cat{X}$, when we invert $\iota_{\cat{K}/}$, are the same.
  For this, it suffices to construct a morphism between these two cospans in $\Cat/\cat{X}$, i.e. a functor $T\colon E/\cat{X}\to D/\cat{X}$ such that
  \[
    \begin{tikzcd}[sep=small]
      & E/\cat{X}
        \ar[dd,"T"]
    \\\cat{K}
        \ar[ur,"\iota_{E/}"]
        \ar[dr,swap,"\widehat{\rho}"]
      && \cat{L}
        \ar[ul,swap,"\widehat{\sigma}"]
        \ar[dl,"\widehat{\sigma\rho}"]
    \\& D/\cat{X}
    \end{tikzcd}
  \]
  commutes. Indeed, from $T \iota_{E/}=\hat \rho$ we will conclude
  $L(\hat \rho) L(\iota_{E/})^{-1}=L(T)$ and thus, from $T \hat \sigma=\widehat{\sigma\rho}$, we'll find $L(\hat\rho) L(\iota_{E/})^{-1} L(\hat\sigma)=L(T) L(\hat \sigma)=L(\widehat{\sigma \rho})$,
  as desired.

  We now define the necessary functor $T$ as induced by the lax natural
  transformation
	 % https://q.uiver.app/#q=WzAsNixbMCwwLCJLIl0sWzEsMSwiWCJdLFswLDIsIlgiXSxbMSwwLCJKIl0sWzIsMSwiWCJdLFsxLDIsIlgiXSxbMCwxLCJFIiwyXSxbMiwxXSxbMCwzLCJSIiwyXSxbMiw1XSxbNSw0XSxbMSw0XSxbMyw0LCJEIl0sWzEyLDYsIlxccmhvIiwyLHsic2hvcnRlbiI6eyJzb3VyY2UiOjIwLCJ0YXJnZXQiOjIwfX1dXQ==
	
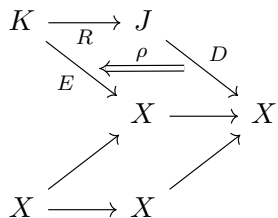
\[\begin{tikzcd}
		K & J \\
		& X & X \\
		X & X
		\arrow[""{name=0, anchor=center, inner sep=0}, "E"', from=1-1, to=2-2]
		\arrow[from=3-1, to=2-2]
		\arrow["R"', from=1-1, to=1-2]
		\arrow[from=3-1, to=3-2]
		\arrow[from=3-2, to=2-3]
		\arrow[from=2-2, to=2-3]
		\arrow[""{name=1, anchor=center, inner sep=0}, "D", from=1-2, to=2-3]
		\arrow["\rho"', shorten <=6pt, shorten >=6pt, Rightarrow, from=1, to=0]
	\end{tikzcd}\]
	of spans, where unlabelled arrows are identities. Explicitly this means we have

  \[
    \pi_{D/} T = R \pi_{E/},\quad \pi_{/\cat{X}} T=\pi_{/\cat{X}},\quad \alpha_{D/\cat{X}}T = \alpha_{E/\cat{X}} \rho,
  \]
  as pictured below.
	% https://q.uiver.app/#q=WzAsNixbMCwwLCJFL1xcY2F0IFgiXSxbMSwwLCJcXGNhdCBYIl0sWzIsMCwiRC9cXGNhdCBYIl0sWzIsMSwiSiJdLFsxLDEsIlxcY2F0IFgiXSxbMCwxLCJLIl0sWzAsMV0sWzAsNSwiXFxwaV97RS99IiwyXSxbMiwzLCJcXHBpX3tEL30iXSxbMiwxXSxbMyw0XSxbNSw0XSxbMSw0LCJcXGlkX3tcXGNhdCBYfSIsMV0sWzAsMiwiVCIsMix7ImN1cnZlIjotM31dLFs1LDMsIlIiLDIseyJjdXJ2ZSI6M31dLFs3LDEyLCJcXGFscGhhX3tFL1xcY2F0IFh9IiwyLHsic2hvcnRlbiI6eyJzb3VyY2UiOjIwLCJ0YXJnZXQiOjMwfX1dLFsxMiw4LCJcXGFscGhhX3tEL1xcY2F0IFh9IiwyLHsic2hvcnRlbiI6eyJzb3VyY2UiOjIwLCJ0YXJnZXQiOjMwfX1dLFsxNCw0LCJcXHJobyIsMix7InNob3J0ZW4iOnsic291cmNlIjoyMH19XV0=
	
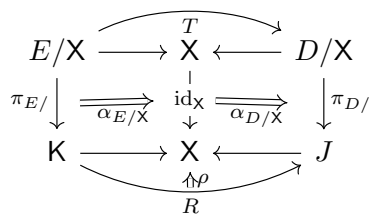
\[\begin{tikzcd}
		{E/\cat{X}} & {\cat{X}} & {D/\cat{X}} \\
		\cat K & {\cat{X}} & J
		\arrow[from=1-1, to=1-2]
		\arrow[""{name=0, anchor=center, inner sep=0}, "{\pi_{E/}}"', from=1-1, to=2-1]
		\arrow[""{name=1, anchor=center, inner sep=0}, "{\pi_{D/}}", from=1-3, to=2-3]
		\arrow[from=1-3, to=1-2]
		\arrow[from=2-3, to=2-2]
		\arrow[from=2-1, to=2-2]
		\arrow[""{name=2, anchor=center, inner sep=0}, "{\id_{\cat{X}}}"{description}, from=1-2, to=2-2]
		\arrow["T"', curve={height=-18pt}, from=1-1, to=1-3]
		\arrow[""{name=3, anchor=center, inner sep=0}, "R"', curve={height=18pt}, from=2-1, to=2-3]
		\arrow["{\alpha_{E/\cat{X}}}"', shorten <=8pt, shorten >=13pt, Rightarrow, from=0, to=2]
		\arrow["{\alpha_{D/\cat{X}}}"', shorten <=9pt, shorten >=13pt, Rightarrow, from=2, to=1]
		\arrow["\rho"', shorten <=1pt, Rightarrow, from=3, to=2-2]
	\end{tikzcd}\]

  To see that $T$ indeed gives a morphism of cospans, note firstly that
  $T \hat \sigma\colon\cat L\to D/\cat{X}$ is induced by the pasting
  % https://q.uiver.app/#q=WzAsNyxbMSwwLCJLIl0sWzIsMSwiWCJdLFsxLDIsIlgiXSxbMiwwLCJKIl0sWzMsMSwiWCJdLFsyLDIsIlgiXSxbMCwxLCJMIl0sWzAsMSwiRSIsMl0sWzIsMV0sWzAsMywiUiIsMl0sWzIsNV0sWzUsNF0sWzEsNF0sWzMsNCwiRCJdLFs2LDAsIlMiLDJdLFs2LDIsIkYiLDJdLFsxMyw3LCJcXHJobyIsMix7InNob3J0ZW4iOnsic291cmNlIjoyMCwidGFyZ2V0IjoyMH19XSxbNywxNSwiXFxzaWdtYSIsMCx7InNob3J0ZW4iOnsic291cmNlIjo0MCwidGFyZ2V0IjoyMH19XV0=
	
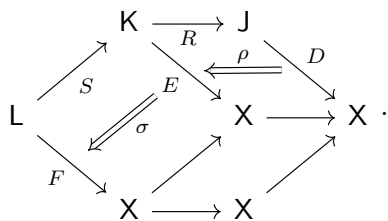
\[\begin{tikzcd}
		& \cat K & \cat J \\
		\cat L && \cat{X} & \cat{X} \\
		& \cat{X} & \cat{X}
		\arrow[""{name=0, anchor=center, inner sep=0}, "E"', from=1-2, to=2-3]
		\arrow[from=3-2, to=2-3]
		\arrow["R"', from=1-2, to=1-3]
		\arrow[from=3-2, to=3-3]
		\arrow[from=3-3, to=2-4]
		\arrow[from=2-3, to=2-4]
		\arrow[""{name=1, anchor=center, inner sep=0}, "D", from=1-3, to=2-4]
		\arrow["S"', from=2-1, to=1-2]
		\arrow[""{name=2, anchor=center, inner sep=0}, "F"', from=2-1, to=3-2]
		\arrow["\rho"', shorten <=6pt, shorten >=6pt, Rightarrow, from=1, to=0]
		\arrow["\sigma", shorten <=14pt, shorten >=7pt, Rightarrow, from=0, to=2]
	\end{tikzcd}.\]
	Thus we have $\pi_{D/} T \hat\sigma = R S,\pi_{/X} T \hat\sigma = F$, and $\alpha_{D/X}* (T\hat \sigma): DRS\To F = \sigma  (\rho  S)$, which is precisely the definition of $\widehat{\sigma\rho}.$

  The other commutative triangle is similar. We see that $T \iota_{E/}$
  is induced from the pasting
  % https://q.uiver.app/#q=WzAsNyxbMSwwLCJLIl0sWzIsMSwiWCJdLFsxLDIsIlgiXSxbMiwwLCJKIl0sWzMsMSwiWCJdLFsyLDIsIlgiXSxbMCwxLCJLIl0sWzAsMSwiRSIsMl0sWzIsMV0sWzAsMywiUiIsMl0sWzIsNV0sWzUsNF0sWzEsNF0sWzMsNCwiRCJdLFs2LDBdLFs2LDIsIkUiLDJdLFsxMyw3LCJcXHJobyIsMix7InNob3J0ZW4iOnsic291cmNlIjoyMCwidGFyZ2V0IjoyMH19XV0=
	
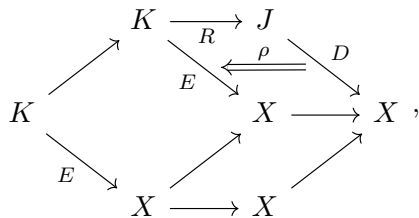
\[\begin{tikzcd}
		& K & J \\
		K && X & X \\
		& X & X
		\arrow[""{name=0, anchor=center, inner sep=0}, "E"', from=1-2, to=2-3]
		\arrow[from=3-2, to=2-3]
		\arrow["R"', from=1-2, to=1-3]
		\arrow[from=3-2, to=3-3]
		\arrow[from=3-3, to=2-4]
		\arrow[from=2-3, to=2-4]
		\arrow[""{name=1, anchor=center, inner sep=0}, "D", from=1-3, to=2-4]
		\arrow[from=2-1, to=1-2]
		\arrow["E"', from=2-1, to=3-2]
		\arrow["\rho"', shorten <=6pt, shorten >=6pt, Rightarrow, from=1, to=0]
	\end{tikzcd},\]
	which also yields $\hat\rho.$ Thus $P$ is indeed a functor.

  \medskip
  \emph{Step 3: $P$ lifts $M$ across $i'$ and extends $L$ along $i$}

  Finally we need to show that the functor $P\colon\DiagOp(\cat{X})\to(\Cat/\cat{X})^\op[(\eosclass^\strict)^{-1}]$ thus defined does indeed make Diagram~\labelcref{equation:diagram-for-P} commute.
  We defined $P$ exactly to get $i' P=M$, so it remains only to check that $Pi=L$.

  If $R\colon(\cat{J},D)\to(\cat{K},E)$ is a morphism in $(\Cat/\cat{X})^\op$, then by definition of $P$ we have
  \[
    \begin{aligned}
      P(i(R))
      &= P(R,\id_E)
    \\&= L(\widehat{\id_E}) L(\iota_{D/})^{-1}.
    \end{aligned}
  \]

  So we want to show that $L(R)=L(\widehat{\id_E}) L(\iota_{D/})^{-1}$, or, equivalently, that
  \[
    L(R) L(\iota_{D/})
    = L(\widehat{\id_E}).
  \]
  But $L(R) L(\iota_{D/})=L(R\iota_{D/})$, and so it suffices to show that (switching to $\Cat$ from $(\Cat/\cat{X})^\op$) we have $\iota_{D/} R=\widehat{\id_E}$ as functors $\cat{K}\to D/\cat{X}$.

  We know that $\widehat{\id_E}$ is determined by $E,R$, and $\id_E.$
  Meanwhile $\iota_{D/} R$ is determined by the following diagram:

  % https://q.uiver.app/#q=WzAsNSxbMCwxLCJLIl0sWzEsMSwiSiJdLFsyLDAsIkoiXSxbMywxLCJYIl0sWzIsMiwiWCJdLFswLDEsIlIiXSxbMSwyXSxbMiwzLCJEIl0sWzEsNCwiRCJdLFs0LDNdXQ==
\[\begin{tikzcd}
	&& J \\
	K & J && X \\
	&& X
	\arrow["R", from=2-1, to=2-2]
	\arrow[from=2-2, to=1-3]
	\arrow["D", from=1-3, to=2-4]
	\arrow["D", from=2-2, to=3-3]
	\arrow[from=3-3, to=2-4].
\end{tikzcd}\]
Since $D R=E$, we see $\widehat{\id_E}=\iota_{D/} R$, as desired.
\end{proof}

\begin{corollary}
\label{corollary:localisation-step-3}
  The canonical functor
  \[
    (\Cat/\cat{X})^\op[(\eosclass^\strict)^{-1}]
    \to \DiagOp(\cat{X})[(\eosclass^\strict)^{-1}]
  \]
  is an isomorphism.
\end{corollary}

\begin{proof}
  First of all, note that this canonical functor is exactly the functor $i'$ in Diagram~\labelcref{equation:diagram-for-P}.
  Since $P$ lands in $(\Cat/\cat{X})^\op[(\eosclass^\strict)^{-1}]$, by definition it sends every morphism in $\eosclass^\strict$ to an isomorphism;
  by the universal property of localisations, it thus factors as $P=P'M$, as in Diagram~\labelcref{equation:diagram-for-P'}.
  It is then purely a matter of diagram chasing and abstract nonsense (which we spell out below) to show that this $P'$ is the inverse to $i'$, witnessing it as an isomorphism.

  We need to show that $P'i'$ and $i'P'$ are both identity functors, but, again by the universal property of localisations, it suffices to show that $P'i'L=L$ and $i'P'M=M$.
  By the commutativity of Diagram~\labelcref{equation:diagram-for-P} and the factorisation $P=P'M$, we see that $P'i'L=P'Mi=Pi$;
  the factorisation $P=P'M$ also implies that $i'P'M=i'P$.
  In summary then, we need to show that
  \[
    \begin{aligned}
      Pi &= L
    \\i'P &= M
    \end{aligned}
  \]
  but both of these follow immediately from the commutativity of Diagram~\labelcref{equation:diagram-for-P} (\cref{lemma:diagram-for-P-commutes}).
\end{proof}

\subsection*{Step 4. Strict weak equivalences and weak equivalences in $\DiagOp(\cat{X})$}

This step is handled rather anticlimactically relative to the previous one.

\begin{corollary}
\label{corollary:localisation-step-4}
  The canonical functor
  \[
    \DiagOp(\cat{X})[(\eosclass^\strict)^{-1}]
    \to \DiagOp(\cat{X})[\eosclass^{-1}]
  \]
  is an isomorphism.
\end{corollary}

\begin{proof}
  Recall the canonical factorisation
  \[
    \begin{tikzcd}
    % https://q.uiver.app/#q=WzAsMyxbMCwwLCJcXGNhdHtKfSJdLFsyLDAsIlxcY2F0e0t9Il0sWzEsMSwiXFxjYXR7WH0iXSxbMSwwLCJSIiwyXSxbMCwyLCJEIiwyXSxbMSwyLCJFIl0sWzQsNSwiXFxyaG8iLDEseyJvZmZzZXQiOi0yLCJzaG9ydGVuIjp7InNvdXJjZSI6MjAsInRhcmdldCI6MjB9fV1d
      {\cat{J}} && {\cat{K}} \\
      & {\cat{X}}
      \arrow["R"', from=1-3, to=1-1]
      \arrow[""{name=0, anchor=center, inner sep=0}, "D"', from=1-1, to=2-2]
      \arrow[""{name=1, anchor=center, inner sep=0}, "E", from=1-3, to=2-2]
      \arrow["\rho"{description}, shift left=2, shorten <=7pt, shorten >=7pt, Rightarrow, from=0, to=1]
    \end{tikzcd}
    \quad=\quad
    \begin{tikzcd}
    % https://q.uiver.app/#q=WzAsNCxbMCwwLCJcXGNhdHtKfSJdLFsyLDAsIlxcY2F0e0t9Il0sWzEsMSwiXFxjYXR7WH0iXSxbMSwwLCJEL1xcY2F0e1h9Il0sWzAsMiwiRCIsMl0sWzEsMiwiRSJdLFszLDAsIlxccGlfe1xcY2F0e0p9L30iLDJdLFsxLDMsIlxcd2lkZWhhdHtcXHJob30iLDJdLFszLDIsIlxccGlfey9cXGNhdHtYfX0iXSxbNCw4LCJcXG9tZWdhIiwwLHsib2Zmc2V0IjotMiwic2hvcnRlbiI6eyJzb3VyY2UiOjIwLCJ0YXJnZXQiOjIwfX1dXQ==
      {\cat{J}} & {D/\cat{X}} & {\cat{K}} \\
      & {\cat{X}}
      \arrow[""{name=0, anchor=center, inner sep=0}, "D"', from=1-1, to=2-2]
      \arrow["E", from=1-3, to=2-2]
      \arrow["{\pi_{D/}}"', from=1-2, to=1-1]
      \arrow["\widehat{\rho}"', from=1-3, to=1-2]
      \arrow[""{name=1, anchor=center, inner sep=0}, "{\pi_{/\cat{X}}}", from=1-2, to=2-2]
      \arrow["\omega", shift left=2, shorten <=4pt, shorten >=4pt, Rightarrow, from=0, to=1]
    \end{tikzcd}
  \]
  from the proof of \cref{lemma:diagram-for-P-commutes}.
  By \cref{corollary:projection-from-comma-is-equivalence}, $\pi_{D/}\colon D/\cat{X}\to\cat{J}$ is a weak equivalence, so 2-out-of-3 for weak equivalences says that, if $(R,\rho)$ is a weak equivalence, then so too is $(\widehat{\rho},\id)$.
  Since $(\widehat{\rho},\id)$ is strict, to show that inverting strict weak equivalences inverts $R$, it suffices to show that inverting strict weak equivalences inverts $(\pi_{D/},\omega)$.
  But this follows immediately from the fact that $\pi_{D/}$ is split by $\iota_{D/}$, which is strict (\cref{lemma:projections-from-comma-are-split-adjoints}) and a weak equivalence (\cref{corollary:projection-from-comma-is-equivalence}).
\end{proof}

\section{The general 2-categorical story}
\label{section:2-categorical-story}

We now recapitulate the argument of \cref{section:localisations-of-categories-of-diagrams} in a more abstract setting.
To start, we recall how completeness and cocompleteness work for $2$-categories.

Let $\twocat{K}$
denote an arbitrary finitely bicomplete $2$-category. We take
limits and colimits in $2$-categories in the enriched sense, so that, for
instance, for a product we require a $2$-natural \emph{isomorphism} of
\emph{categories}
$\twocat{K}(A,B\times C)\cong \twocat{K}(A,B)\times \twocat{K}(A,C)$. Thus a pair of
$1$-morphisms into $B$ and $C$ lift fully uniquely to a $1$-morphism into $B\times C$,
and a $2$-morphism between $f,g\colon A\to B\times C$ is uniquely determined by its
whiskerings with the projections to $B$ and $C$.

Finite completeness, in the enriched sense, amounts to the assumption that $\twocat{K}$ has finite products and equalisers as well as \emph{cotensors} by the interval
category $\mathbb{2}=0\to 1$.
The cotensor $b^\mathbb{2}$ is a representation of the
arrow category of a hom-category:
$\twocat{K}(A,B^\mathbb{2})\cong \twocat{K}(A,B)^\mathbb{2}$. Thus a $1$-morphism
$\alpha\colon A\to B^\mathbb{2}$ corresponds uniquely to a $2$-morphism between 1-morphisms from $A$ to $B$, which,
by abuse of notation, we'll also denote by $\alpha\colon f\To g:A\to B$. Then a $2$-morphism
$\alpha\To\beta:A\to B^\mathbb{2}$ corresponds to a commutative square
in $\twocat{K}(A,B)$ as below:
\[
% https://q.uiver.app/#q=WzAsNCxbMCwwLCJmIl0sWzEsMCwiZyJdLFswLDEsImgiXSxbMSwxLCJrIl0sWzAsMSwiXFxhbHBoYSJdLFsyLDMsIlxcYmV0YSJdLFswLDJdLFsxLDNdXQ==
  \begin{tikzcd}
  	f & g \\
  	h & k
  	\arrow["\alpha", from=1-1, to=1-2]
  	\arrow["\beta", from=2-1, to=2-2]
  	\arrow[from=1-1, to=2-1]
  	\arrow[from=1-2, to=2-2]
  \end{tikzcd}
\]

Dually, $\twocat{K}$ is finitely cocomplete when it has finite coproducts,
coequalisers, and \emph{tensors} by $\mathbb{2}$, which by definition
satisfy $\twocat{K}(A\otimes\mathbb{2},B)\cong\twocat{K}(A,B)^{\mathbb{2}}$.

\begin{example}\label{example:2-categories-of-interest}
  Many $2$-categories of interest are finitely bicomplete.
  For example:

  \begin{itemize}
    \item $\Cat$ is finitely bicomplete, with tensors given
    by the cartesian product and cotensors by the internal hom.
    \item The $2$-categories of multicategories, symmetric multicategories,
      cartesian multicategories, polycategories, etc. are also finitely bicomplete. Indeed, all these categories are finitely locally presentable, since they are the categories of models of essentially algebraic theories.

      For instance, if $X$ is a multicategory,
      then the cotensor $X^{\mathbb{2}}$ is the multicategory whose objects are the \emph{unary}
      morphisms of $X$ and whose multimorphisms $X^{\mathbb{2}}((f),g)$ are commutative squares of the form
      \[
        \begin{tikzcd}
          (x)
            \ar[r]
            \ar[d,swap,"(f)"]
          & z
            \ar[d]
        \\(y)
            \ar[r]
          & w
        \end{tikzcd}.
      \]
      It is interesting to note that the $2$-category of multicategories is \emph{not} monoidal closed in any obvious way specializing to its cotensoring by
      categories.
    \item The $2$-category of strict algebras for a reasonable $2$-monad on a reasonable $2$-category (say, an accessible $2$-monad on a locally presentable $2$-category, see \cite[Theorem 3.8]{blackwell1989})
      and \emph{strict} morphisms is finitely bicomplete.
      For instance, this includes the $2$-category of monoidal categories and
      strict monoidal functors.
  \end{itemize}

  We are, naturally, more interested in $2$-categories like that of
  monoidal categories and strong monoidal functors, cartesian categories
  and functors preserving finite products in the usual sense,
  multicategories with finite products, and other
  examples of 2-categories of algebras and \emph{pseudo} morphisms over a 2-monad.
  Such 2-categories generally possess only \emph{flexible} limits and
  \emph{shrinkable} colimits, in the sense of \cite{bourke2023}.

  For instance, there is no initial object
  in the $2$-category of monoidal categories, nor any equaliser of the two
  maps of cartesian categories from the terminal category to the isomorphism
  category $\mathbb{I}$. We shall later handle these cases by showing
  that we can cover our intended applications while considering only
  the strict morphisms.
\end{example}

We now define a notion of discrete opfibration in an arbitrary $2$-category that is generally quite useful.

\begin{definition}
\label{definition:2-categorical-dopf}
  A morphism $p\colon E\to B$ in a $2$-category $\twocat{K}$ is called a
  \emph{discrete opfibration} if, given a $2$-morphism
  $\alpha\colon f\To g:X\to B$ and a $1$-morphism $\overline{f}\colon X\to E$ with
  $\overline{f}\cdot p=f$, there exists a unique $2$-morphism
  $\overline{\alpha}\colon\overline{f}\To\overline{g}:X\to E$ with $\overline{\alpha}\cdot p=\alpha$.
\end{definition}

If $\twocat{K}$ admits tensors with $\mathbb{2}$, then this definition
can be re-expressed as the existence, given the solid square below,
of a unique $1$-morphism with the signature of the dotted arrow making
both triangles commute:
\begin{equation*}
  % https://q.uiver.app/#q=WzAsNCxbMCwwLCJYIl0sWzEsMCwiRSJdLFswLDEsIlhcXG90aW1lc1xcbWF0aGJiezJ9Il0sWzEsMSwiQiJdLFswLDEsIlxcYmFyIGYiXSxbMCwyLCJpXzAiLDJdLFsxLDMsInAiXSxbMiwzXSxbMiwxLCIiLDAseyJzdHlsZSI6eyJib2R5Ijp7Im5hbWUiOiJkYXNoZWQifX19XV0=
  \begin{tikzcd}
    X & E \\
    {X\otimes\mathbb{2}} & B
    \arrow["{\bar f}", from=1-1, to=1-2]
    \arrow["{i_0}"', from=1-1, to=2-1]
    \arrow["p", from=1-2, to=2-2]
    \arrow[from=2-1, to=2-2]
    \arrow[dashed, from=2-1, to=1-2]
  \end{tikzcd}
\end{equation*}
where the map $i_0\colon X\to X\otimes\mathbb{2}$ is
induced by the inclusion $0:\mathbb 1\to \mathbb 2$ of the initial object into $\mathbb 2.$

In the language of factorisation systems, discrete opfibrations are thus
precisely the morphisms orthogonal on the right to every morphism
$i_0\colon X\to X\otimes\mathbb{2}$.
A few remarks on this definition are in order.

\begin{remark}
\label{remark:about-2-categorical-dopfs}
  \begin{enumerate}[i.]
    \item It is easy to check that a discrete opfibration of categories in the
      ordinary sense (\cref{definition:discrete-opfibration}) induces a discrete opfibration in the abstract sense above (\cref{definition:2-categorical-dopf}).
    \item We do not claim that \cref{definition:2-categorical-dopf} is the best notion of discrete opfibration
      for \emph{all} purposes and \emph{all} $2$-categories $\twocat{K}$. A key 2-category of interest to us
      is that of multicategories, and there there exists a natural stricter notion of
      discrete opfibration, allowing for lifts against \emph{multi}morphisms
      with a lifted domain. However, the notion used here suffices to prove our main theorem for multicategories.
      We thus will be able to describe the localization of the category of diagrams in a multicategory at weak
      equivalences defined with respect to this broader notion of discrete opfibration; we leave open the question
      of whether there are interesting examples of diagrams in a multicategory only weakly equivalent with respect to the stricter
      discrete opfibrations.
    \item While this notion of discrete opfibration for multicategories captures \emph{at least}
      everything we intend to capture, for other $2$-categories it is well-known that discrete
      opfibrations capture very \emph{little} of interest. If $\twocat{K}=\enrichedcat{\cat{V}}$ for a
      non-cartesian symmetric monoidal $\cat{V}$, then, in general, discrete opfibrations do not fully
      capture $\cat{V}$-functors into $\cat{V}$, which are instead modelled using discrete
      op-\emph{co}fibrations, and therefore discrete opfibrations and thus the result to be proven are of limited interest;
      in short, the story of a diagram in such a $\cat V$-category as a system of equations to be solved by lifting against
      a discrete opfibration does not go through.
  \end{enumerate}
  \,
\end{remark}

Let us exhibit the key general example of discrete opfibrations.

\begin{proposition}
\label{proposition:coslice-is-dopf}
  Given any cospan $1 \xrightarrow{x} A \xleftarrow{f} B$ in a finitely complete $2$-category
  $\twocat{K}$, the projection from the comma object $x/ f\to B$
  is a discrete opfibration.
\end{proposition}

\begin{proof}
  Recall that the comma object is the terminal inhabitant of the right-hand side of the following situation:
% https://q.uiver.app/#q=WzAsNixbMSwwLCJ4XFxkb3duYXJyb3cgZiJdLFsyLDAsIjEiXSxbMSwxLCJBIl0sWzIsMSwiQiJdLFswLDAsIkMiXSxbMCwxLCJDXFxvdGltZXNcXG1hdGhiYnsyfSJdLFswLDFdLFswLDJdLFsxLDMsIngiXSxbMiwzLCJnIiwyXSxbMSwyLCIiLDIseyJzaG9ydGVuIjp7InNvdXJjZSI6MzAsInRhcmdldCI6MzB9LCJsZXZlbCI6Mn1dLFs0LDBdLFs1LDJdLFs0LDVdXQ==
\[\begin{tikzcd}
	C & {x/ f} & 1 \\
	{C\otimes\mathbb{2}} & A & B
	\arrow[from=1-2, to=1-3]
	\arrow[from=1-2, to=2-2]
	\arrow["x", from=1-3, to=2-3]
	\arrow["g"', from=2-2, to=2-3]
	\arrow[shorten <=6pt, shorten >=6pt, Rightarrow, from=1-3, to=2-2]
	\arrow[from=1-1, to=1-2]
	\arrow[from=2-1, to=2-2]
	\arrow[from=1-1, to=2-1]
\end{tikzcd}\]

Given maps from $C$ and $C\otimes\mathbb 2$ as shown, we must show there exists a unique lift $C\otimes \mathbb 2\to x/ f.$
The given information amounts to a 2-morphism $f:f_0\To f_1:C\to A$ and a 2-morphism $t:!_C\cdot x\To f_0\cdot g:C\to B$,
where $!_C:C\to 1$ is the unique map to the terminal object. The desired map will be uniquely determined, by the universal property of
$x/ f$, by a 2-morphism $!_{C\otimes\mathbb 2}\cdot x\To f\cdot g: C\otimes 2\to B.$ Such a 2-morphism is itself given,
by the universal property of the tensor product, by a choice of 2-morphism $u:!_C\cdot x\To f_1\cdot g$
such that the square below commutes. But such a choice is visibly unique since the top leg of the square is an identity.
% https://q.uiver.app/#q=WzAsNCxbMCwwLCJ4XFxjaXJjICFfQyJdLFsxLDAsInhcXGNpcmMgIV9DIl0sWzAsMSwiZ1xcY2lyYyBmXzAiXSxbMSwxLCJnXFxjaXJjIGZfMSJdLFswLDEsIlxcaWQiXSxbMCwyLCJ0IiwyXSxbMiwzLCJnXFxjaXJjIGYiLDJdLFsxLDMsInUiLDAseyJzdHlsZSI6eyJib2R5Ijp7Im5hbWUiOiJkb3R0ZWQifX19XV0=
\[\begin{tikzcd}
	{!_C\cdot x} & {!_C\cdot x} \\
	{ g\cdot f_0} & {f_1\cdot g}
	\arrow["\id", from=1-1, to=1-2]
	\arrow["t"', from=1-1, to=2-1]
	\arrow["{ g\cdot g}"', from=2-1, to=2-2]
	\arrow["u", dotted, from=1-2, to=2-2]
\end{tikzcd}\]

\end{proof}

This result generalizes the analogous fact for coslice categories, where $f=\id_B.$
However, the sense in which $x/ f$ generalizes a coslice category can
be slightly unintuitive, so let us explain it in more detail in the case of multicategories.
The terminal multicategory $1$ has a single object $\star$ and a single $n$-ary
morphism $(\star)\to \star$ for every $n$.
So if $X$ is a multicategory, a
multifunctor $x\colon 1\to X$ picks out a \emph{monoid} in $X$, not just an object.

If $m$ is a monoid in $X$, then the objects of the coslice multicategory $m/X$ are unary
morphisms $m\to a$, and its multimorphisms $(m/X)((m\to a),m\to b)$ are given by
morphisms $(a)\to b$ such that a square involving the monoid operation of
$m$ commutes. In contrast, if $a$ is a mere object of a multicategory, it
is represented by a functor from the multicategory $1_\bot$ with one object
and no non-unary morphisms, and then the analogous comma multicategory
$a/X$ also has no non-unary morphisms.

A particularly important example for applications is when $X$ is a multicategory
of vector spaces, vector bundles, or sheaves of vector spaces on some manifold.
In any such case the
standard notion of solution of a (system of) equation(s) corresponds to lifting
against the codomain projection out of $\mathbb{R}/X$, where $\mathbb{R}$ represents
respectively the one-dimensional real vector space, the rank-one trivial bundle,
or the sheaf of continuous real-valued functions. In every case $\mathbb{R}$ is a
monoid and so the full multicategorical slice exists, allowing us to continue
with the story of solving equations in terms of lifting against
discrete opfibrations.

From any notion of discrete opfibration, we immediately get a corresponding notion of initiality:

\begin{definition}
\label{definition:initial-morphism}
  A morphism $f\colon X\to Y$ in a $2$-category $\twocat{K}$
  is called \emph{initial} if it is left orthogonal to all discrete opfibrations:
  thus if, given the solid square below in which $p$ is a discrete opfibration,
  there exists a unique diagonal lift as indicated, making both triangles commute:
  \[
    \begin{tikzcd}
      X
        \ar[r]
        \ar[d,swap,"f"]
      & E
        \ar[d,"p"]
    \\Y
        \ar[r]
        \ar[ur,dashed]
      & B
    \end{tikzcd}
  \]
\end{definition}

We are aiming to construct a comprehensive factorisation system on $\twocat{K}$,
and this property must hold for any putative class of initial maps --- to ensure that this is the case, we have simply baked it in to the definition.
While this is less satisfying than the independent
definition for initial functors between categories,
we can nevertheless show that initial $1$-morphisms have a familiar sufficient
classification:

\begin{proposition}
\label{proposition:left-adjoint-is-initial}
  Every left adjoint in a $2$-category $\twocat{K}$ is initial.
\end{proposition}

\begin{proof}
  Let $\ell\colon X\to Y$ be a left adjoint in $\twocat{K}$, so we can choose
  $r\colon Y\to X$, $\eta\colon 1_X\To\ell\cdot r$, and $\varepsilon\colon r\cdot\ell\To 1_Y$ satisfying
  the usual triangle equations. Consider also a discrete opfibration $p:E\to B$
  in $\twocat{K}$, and a commutative square
  % https://q.uiver.app/#q=WzAsNCxbMCwwLCJYIl0sWzAsMSwiWSJdLFsxLDAsIkUiXSxbMSwxLCJCIl0sWzAsMSwiXFxlbGwiLDJdLFsyLDMsInAiXSxbMCwyLCJmIl0sWzEsMywiZyIsMl1d
  \[\begin{tikzcd}
  	X & E \\
  	Y & B
  	\arrow["\ell"', from=1-1, to=2-1]
  	\arrow["p", from=1-2, to=2-2]
  	\arrow["f", from=1-1, to=1-2]
  	\arrow["g"', from=2-1, to=2-2]
  \end{tikzcd}\]
  We are going to construct unique a $1$-morphism $Y\to E$ factoring $\ell$ through $f$, and factoring $g$
  through $p$.

  \medskip

  \emph{Uniqueness.}
  Suppose we have two diagonal fillers $k,k'$, so that
  % https://q.uiver.app/#q=WzAsNCxbMCwwLCJYIl0sWzAsMSwiWSJdLFsxLDAsIkUiXSxbMSwxLCJCIl0sWzAsMSwiXFxlbGwiLDJdLFsyLDMsInAiXSxbMCwyLCJmIl0sWzEsMywiZyIsMl0sWzEsMiwiayciLDEseyJvZmZzZXQiOjJ9XSxbMSwyLCJrIiwxLHsib2Zmc2V0IjotMn1dXQ==
  \[
    \begin{tikzcd}
    	X & E \\
    	Y & B
    	\arrow["\ell"', from=1-1, to=2-1]
    	\arrow["p", from=1-2, to=2-2]
    	\arrow["f", from=1-1, to=1-2]
    	\arrow["g"', from=2-1, to=2-2]
    	\arrow["{k'}", swap, shift right=1, from=2-1, to=1-2]
    	\arrow["k", shift left=1, from=2-1, to=1-2]
    \end{tikzcd}
  \]
  In particular then, we have $\ell k=\ell k'=f$, and thus
  $r\ell k=r\ell k'=rf$. Since $kp=k'p=g$, we have
  $\varepsilon k p=\varepsilon k' p=\varepsilon g\colon rlg \To g$.
  Since $p$ is a discrete opfibration, there can exist only one
  $2$-morphism $\alpha$ of domain $rf$ such that $\alpha p = \varepsilon g$.
  We have two candidates $\varepsilon k: rf \To k$ and $\varepsilon k': rf \To k'$ for
  such a $2$-morphism,
  which implies that $\varepsilon k=\varepsilon k'$. In particular,
  the codomains of these $2$-morphisms must be equal, whence $k=k'$.

  \medskip

  \emph{Existence.}
  There is a $1$-morphism $Y\to E$ ready to hand, namely,
  $r\cdot f$. Unfortunately, this is not quite sufficient, since
  the composite $\ell r f$ need not coincide with $f$, nor
  $rfp=r\ell g$ with $g$. Noting the irresistible fact
  that $rf$ is correct ``modulo'' the unit and counit
  of the adjunction, however, we look to deform
  $rf$ into the desired filler.

  To construct the filler, we apply the discrete opfibration property to the following situation:
  % https://q.uiver.app/#q=WzAsNCxbMCwxLCJZIl0sWzEsMSwiWSJdLFsyLDEsIkIiXSxbMiwwLCJFIl0sWzAsMSwiXFxpZF9ZIiwyLHsiY3VydmUiOjJ9XSxbMCwxLCJyXFxlbGwiLDAseyJjdXJ2ZSI6LTJ9XSxbMSwyLCJnIiwyXSxbMywyLCJwIl0sWzAsMywicmYiLDEseyJjdXJ2ZSI6LTN9XSxbNSw0LCJcXHZhcmVwc2lsb24iLDIseyJzaG9ydGVuIjp7InNvdXJjZSI6MjAsInRhcmdldCI6MjB9fV1d
  \[
    \begin{tikzcd}
    	&& E \\
    	Y & Y & B
    	\arrow[""{name=0, anchor=center, inner sep=0}, "{\id_Y}"', curve={height=12pt}, from=2-1, to=2-2]
    	\arrow[""{name=1, anchor=center, inner sep=0}, "r\ell", curve={height=-12pt}, from=2-1, to=2-2]
    	\arrow["g"', from=2-2, to=2-3]
    	\arrow["p", from=1-3, to=2-3]
    	\arrow["rf", curve={height=-18pt}, from=2-1, to=1-3]
    	\arrow["\varepsilon"', shorten <=3pt, shorten >=3pt, Rightarrow, from=1, to=0]
    \end{tikzcd}
  \]
This produces a 2-morphism
$\alpha\colon rf\To k:Y\to E$ such that $\alpha p=\varepsilon g$.
In particular, the codomain $kp$ of $\alpha p$ must
coincide with the codomain $g$ of $\varepsilon g$, and so
we see that $k$ already satisfies the desideratum $kp=g.$

  Now consider the following diagram:
% https://q.uiver.app/#q=WzAsNCxbMCwwLCJYIl0sWzIsMCwiRSJdLFsyLDIsIkIiXSxbMCwyLCJYIl0sWzAsMSwiZiIsMSx7ImN1cnZlIjotNH1dLFswLDEsIlxcZWxsIGsiLDEseyJjdXJ2ZSI6NH1dLFswLDEsIlxcZWxsIHIgZiIsMV0sWzEsMiwicCJdLFswLDMsIlxcaWRfWCIsMl0sWzMsMiwiXFxlbGwgZyIsMix7ImN1cnZlIjo0fV0sWzMsMiwiXFxlbGwgZyIsMSx7ImN1cnZlIjotNH1dLFszLDIsIlxcZWxsIHJcXGVsbCBnIiwxXSxbNCw2LCJcXGV0YSBmIiwwLHsic2hvcnRlbiI6eyJzb3VyY2UiOjIwLCJ0YXJnZXQiOjIwfX1dLFs2LDUsIlxcZWxsXFxhbHBoYSIsMCx7InNob3J0ZW4iOnsic291cmNlIjoyMCwidGFyZ2V0IjoyMH19XSxbMTAsMTEsIlxcZXRhIFxcZWxsIGciLDIseyJzaG9ydGVuIjp7InNvdXJjZSI6MjAsInRhcmdldCI6MjB9fV0sWzExLDksIlxcZWxsIFxcdmFyZXBzaWxvbiBnIiwyLHsic2hvcnRlbiI6eyJzb3VyY2UiOjIwLCJ0YXJnZXQiOjIwfX1dXQ==
\[\begin{tikzcd}
	X && E \\
	\\
	X && B
	\arrow[""{name=0, anchor=center, inner sep=0}, "f"{description}, curve={height=-24pt}, from=1-1, to=1-3]
	\arrow[""{name=1, anchor=center, inner sep=0}, "{\ell k}"{description}, curve={height=24pt}, from=1-1, to=1-3]
	\arrow[""{name=2, anchor=center, inner sep=0}, "{\ell r f}"{description}, from=1-1, to=1-3]
	\arrow["p", from=1-3, to=3-3]
	\arrow["{\id_X}"', from=1-1, to=3-1]
	\arrow[""{name=3, anchor=center, inner sep=0}, "{\ell g}"', curve={height=24pt}, from=3-1, to=3-3]
	\arrow[""{name=4, anchor=center, inner sep=0}, "{\ell g}"{description}, curve={height=-24pt}, from=3-1, to=3-3]
	\arrow[""{name=5, anchor=center, inner sep=0}, "{\ell r\ell g}"{description}, from=3-1, to=3-3]
	\arrow["{\eta f}", shorten <=3pt, shorten >=3pt, Rightarrow, from=0, to=2]
	\arrow["\ell\alpha", shorten <=3pt, shorten >=3pt, Rightarrow, from=2, to=1]
	\arrow["{\eta \ell g}"', shorten <=3pt, shorten >=3pt, Rightarrow, from=4, to=5]
	\arrow["{\ell \varepsilon g}"', shorten <=3pt, shorten >=3pt, Rightarrow, from=5, to=3]
\end{tikzcd}\]
By hypothesis, $\eta f p = \eta \ell g$ and $\ell \alpha p=\ell \varepsilon g$. That is,
the diagram commutes serially.

By the triangle identities,
the vertical composite $(\eta \ell g)\cdot (\ell \varepsilon g)$ is equal to $\id_{\ell g}$. In
other words, the vertical composite in the upper half of the diagram is a lift
of $\id_{\ell g}=\id_{fp}$ along $p$ with domain $f$.
But $\id_f$ itself is another such lift, and so,
by the discrete opfibration property, we conclude that
$\id_f=(\eta f)\cdot (\ell\alpha)$. In particular, this
implies that $\ell k=f$, and since we already had $kp=g$, we have established
that $k$ gives the desired diagonal lifting.

\end{proof}

If $X$ is an object of any $2$-category $\twocat{K}$, we can define the categories $\Diag(X)$ and $\DiagOp(X)$ of
diagrams in $X$, as well as their pseudo and strict variants, just as we
did for $X\in\Cat$.
We highlight that $\Diag(X)$ is merely a category,
not an object of $\twocat{K}$. The formal properties of these categories are quite
similar in this generality to what we saw in \cref{section:diagram-categories}. For instance:

\begin{lemma}
  For any discrete opfibration $p\colon E\to X$ in a finitely bicomplete $2$-category $\twocat{K}$,
  there is an induced discrete opfibration $\DiagOp(p)\colon\DiagOp(E)\to\DiagOp(X)$ of categories.
\end{lemma}

\begin{proof}
  The functor itself is defined simply by
  whiskering, so we have only to show it is a
  discrete opfibration.

  As in the case that $\twocat{K}=\Cat$, let $d\colon J\to X$ and $d'\colon J\to X$ be diagrams in
  $X$, and let $(r,\rho)\colon d\to d'$ be a morphism in $\DiagOp(X)$.
  Let $\overline{d}$ be a lift of $d$ along $p$, so that we have the diagram
  % https://q.uiver.app/#q=WzAsNCxbMCwwLCJKIl0sWzIsMCwiSiciXSxbMSwxLCJYIl0sWzEsMiwiRSJdLFsxLDAsInIiLDJdLFswLDIsImQiLDJdLFsxLDIsImQnIl0sWzAsMywiXFxiYXIgZCIsMix7ImN1cnZlIjoyfV0sWzMsMiwicCIsMl0sWzUsNiwiXFxyaG8iLDAseyJzaG9ydGVuIjp7InNvdXJjZSI6MjAsInRhcmdldCI6MjB9fV1d
  \[\begin{tikzcd}
  	J && {J'} \\
  	& X \\
  	& E
  	\arrow["r"', from=1-3, to=1-1]
  	\arrow[""{name=0, anchor=center, inner sep=0}, "d"', from=1-1, to=2-2]
  	\arrow[""{name=1, anchor=center, inner sep=0}, "{d'}", from=1-3, to=2-2]
  	\arrow["{\overline{d}}"', curve={height=12pt}, from=1-1, to=3-2]
  	\arrow["p"', from=3-2, to=2-2]
  	\arrow["\rho", shorten <=6pt, shorten >=6pt, Rightarrow, from=0, to=1]
  \end{tikzcd}\]
  We need to show that there exists a unique lift $\overline{\rho}\colon r\overline{d}\to\overline{d'}$ of $\rho$.

  But we have given a $2$-morphism $\rho:rd\To d'$
  into $X$ together with a lift $r\overline{d}$ of its
  domain, so the unique existence of $\overline{\rho}$
  is a direct application of the definition
  of discrete opfibrations in $\twocat{K}$.
\end{proof}

Given this result, we can define weak
equivalences in $X\in\twocat{K}$ just
as we did in $X\in\Cat$.

\begin{definition}
  A morphism $(r,\rho)\colon d\to d'$ in $\DiagOp(X)$ is a \emph{weak equivalence} if the function $\overline{d}\mapsto\overline{d'}$ constructed above is a bijection for every discrete opfibration $p$,
  and similarly for pseudo morphisms in $\Diag(X).$
\end{definition}

As before, it is straightforward to establish
that pseudo morphisms in $\DiagOp(X)$ are weak
equivalences if and only if
their mates are also weak equivalences in $\Diag(X)$.

The diagram categories also continue to enjoy
the following property:

\begin{lemma}
  Let $(r,\rho)\colon(J,d)\to(J',d')$ be a morphism in $\Diag(X)$ such that $r\colon J\to J'$ is initial.
  Then $(r,\rho)$ is a weak equivalence.
\end{lemma}

\begin{proof}
  Given a lift $\overline{d}$ of $d$ along some discrete opfibration $p\colon E\to X$, we need to show that there exists a unique lift
  $\overline{d'}$ of $d'$ along $p$.

  First, we lift $\rho$ itself. Indeed, since
  we're given a $2$-morphism $\rho\colon d\To rd'$
  into $X$ and a lift $\overline{d}$ of its domain, there is
  a unique $\overline{\rho}\colon \overline{d}\To s$ over
  $\rho$.
  \begin{equation*}
    % https://q.uiver.app/#q=WzAsNCxbMCwwLCJKIl0sWzIsMCwiSiciXSxbMiwxLCJYIl0sWzAsMSwiRSJdLFswLDEsInIiXSxbMCwyLCJkIiwxXSxbMSwyLCJkJyJdLFswLDMsIlxcYmFyIGQiLDIseyJjdXJ2ZSI6Mn1dLFszLDIsInAiLDJdLFswLDMsInMiLDAseyJjdXJ2ZSI6LTJ9XSxbNSw2LCJcXHJobyIsMCx7InNob3J0ZW4iOnsic291cmNlIjoyMCwidGFyZ2V0IjoyMH19XSxbNyw5LCJcXGJhclxccmhvIiwwLHsic2hvcnRlbiI6eyJzb3VyY2UiOjIwLCJ0YXJnZXQiOjIwfX1dXQ==
    \begin{tikzcd}
      J && {J'} \\
      E && X
      \arrow["r", from=1-1, to=1-3]
      \arrow[""{name=0, anchor=center, inner sep=0}, "d"{description}, from=1-1, to=2-3]
      \arrow[""{name=1, anchor=center, inner sep=0}, "{d'}", from=1-3, to=2-3]
      \arrow[""{name=2, anchor=center, inner sep=0}, "{\bar d}"', curve={height=12pt}, from=1-1, to=2-1]
      \arrow["p"', from=2-1, to=2-3]
      \arrow[""{name=3, anchor=center, inner sep=0}, "s", curve={height=-12pt}, from=1-1, to=2-1]
      \arrow["\rho", shorten <=6pt, shorten >=6pt, Rightarrow, from=0, to=1]
      \arrow["\bar\rho", shorten <=5pt, shorten >=5pt, Rightarrow, from=2, to=3]
    \end{tikzcd}
  \end{equation*}

  But now we have a square with edges
  $r$, $d'$, $p$, and $s$, and the definition of initiality
  says precisely that there is a unique
  factorisation $\overline{d'}: J' \to E$ of $s$ through $r$
  which lifts $d'$ \mbox{along $p$}.
\end{proof}

As a final preparation for the main result of this section, we explain sufficient
conditions under which
$\twocat{K}$ admits a comprehensive factorisation
system. Since initial morphisms are left orthogonal to discrete opfibrations
by definition, any factorization of a morphism into an initial followed by
a discrete opfibration will be unique up to unique isomorphism; but to show
that such factorizations exist, we require some stronger assumptions
on $\twocat K$, namely, local presentability.

Recall that local presentability of a $2$-category
is closely related, but not equivalent, to local
presentability of its underlying $1$-category. To wit,
a $2$-category $\twocat{K}$ is locally presentable as a
$2$-category when it is cocomplete and admits a small set
$\cat{G}$ of objects which jointly detects isomorphisms
(not just equivalences) and which are finitely presentable
in the 2-categorical sense. This means that
the representable 2-functor corresponding
to any object $G\in \cat{G}$ sends filtered colimits
in $\twocat{K}$ to filtered colimits in $\Cat$. Since the
forgetful functor $\Cat\to\Set$ preserves filtered
colimits, a locally presentable 2-category is locally
presentable as a 1-category. The converse holds in the
context of 2-categories already assumed cocomplete, but
not in general.

In any event, all of our $2$-categories of interest,
amounting to the algebras and strict morphisms for
accessible $2$-monads on categories, multicategories, or
symmetric multicategories, are locally presentable, and
indeed every specific example of interest is locally
\emph{finitely} presentable.

\begin{lemma}
\label{lemma:abstract-comprehensive-factorisation}
  If $\twocat{K}$ is locally presentable, then
  it admits a factorisation system
  whose left class is the class of initial
  morphisms and whose right class is the
  class of discrete opfibrations.
\end{lemma}

\begin{proof}
  This argument is well-known to experts, particularly
  in the unenriched case, but we give the argument for convenience.
  Let $\twocat{K}$ be locally
  $\lambda$-presentable. The class of discrete opfibrations,
  having been defined in terms of a right orthogonality
  property, is closed under all limits in the arrow
  $2$-category $\twocat{K}^{\mathbb{2}}$. Furthermore, any
  morphism $i_0\colon J\to J\otimes \mathbb{2}$ may
  be written as a $\lambda$-filtered colimit of the $i_0$s
  corresponding to $\lambda$-presentable objects.

  This implies that a morphism $p:E\to B$ in $\twocat{K}$ is a
  discrete opfibration if and only if it has the
  discrete opfibration property with respect to the morphisms
  $i_0\colon J\to J\otimes \mathbb{2}$ such that $J$ is
  $\lambda$-presentable. Indeed, the class of morphisms
  left orthogonal to $p$ is closed under \emph{all} colimits
  in the arrow category, not just $\lambda$-filtered ones.

  Now we can see that discrete opfibrations
  are closed under $\lambda$-filtered colimits in $\twocat{K}^{\mathbb{2}}$, since, given a square
  % https://q.uiver.app/#q=WzAsNCxbMCwwLCJKIl0sWzAsMSwiSlxcb3RpbWVzXFxtYXRoYmIgMiJdLFsxLDAsIkUiXSxbMSwxLCJCIl0sWzAsMSwiaV8wIl0sWzAsMl0sWzIsMywicCIsMl0sWzEsM11d
  \[
    \begin{tikzcd}
    	J & E \\
    	{J\otimes\mathbb{2}} & B
    	\arrow["{i_0}", from=1-1, to=2-1]
    	\arrow[from=1-1, to=1-2]
    	\arrow["p"', from=1-2, to=2-2]
    	\arrow[from=2-1, to=2-2]
    \end{tikzcd}
  \]
  with $J$ $\lambda$-presentable and $p$ a $\lambda$-filtered
  colimit of discrete opfibrations, the whole square can
  be factored through one of the members of the
  $\lambda$-filtered colimit.

  Since the class $\cat{M}$ of discrete opfibrations is
  closed in $\twocat{K}^{\mathbb{2}}$ under limits and
  $\lambda$-filtered colimits, we conclude that $\cat{M}$
  is a reflective subcategory of $\twocat{K}^{\mathbb{2}}$.
  It is then straightforward to check that the unit of the
  reflection gives the desired factorisation.
\end{proof}

We are now in a position to reiterate the proof of
\cref{theorem:localisation-of-slice-cat-and-of-diag-op} in the current more general
situation.

\begin{theorem}
\label{theorem:2-categorical-main-theorem}
  Given a locally presentable $2$-category $\twocat K$ containing an object $X$, let $\iclass$ be the class of morphisms in $\twocat{K}/X$ determined by an initial morphism in $\twocat K.$ If $\eosclass$ denotes the weak equivalences in $\DiagOp(X)$, then the canonical functor $(\twocat{K}/X)[\iclass^{-1}]\to \DiagOp(X)[\eosclass^{-1}]$ is an isomorphism.
\end{theorem}

\begin{proof}
  We follow the same four steps as in \cref{section:localisations-of-categories-of-diagrams} for the proof of \cref{theorem:localisation-of-slice-cat-and-of-diag-op}.

  Step 1 really applies to any factorisation
  system in any category.
  In particular, since $\twocat{K}$ is locally presentable, there is a comprehensive
  factorisation system, and so Step 1 goes through as desired.

  Step 2 is handled as in \cref{section:localisations-of-categories-of-diagrams}.

  For Step 3, we first note that the properties given in \cref{lemma:slice-into-diag-has-adjoint} follow just from the universal property
  of the comma object, and not on its concrete description in $\Cat$. The same
  is true for \cref{lemma:projections-from-comma-are-split-adjoints}. Once
  more, the proof of \cref{lemma:diagram-for-P-commutes} makes no use
  of the fact that, there, we had $\twocat{K}=\Cat$.

  For Step 4, again, there is nothing new to say.
\end{proof}

To close, we explain why the above result, which computes
the localisations of the diagram categories in a
locally presentable $2$-category, actually can be extended to handle every
$2$-categories like that of categories with chosen finite
products and functors preserving these up to isomorphism,
which is by no means locally presentable.

In all these examples, we are interested in the $2$-category
of algebras for an accessible $2$-monad $T$ on a locally
presentable $2$-category $\twocat{K}$ (namely, $\twocat{K}$ is
either $\Cat$ or some flavor of multicategory.) As is
well-known, in this case, the $2$-category $\algebra{T}_\strict$ of $T$-algebras
and \emph{strict} morphisms
is again locally presentable. Thus \cref{theorem:2-categorical-main-theorem} applies perfectly well to $\algebra{T}_\strict$.

Really, we are interested in $\algebra{T}$ itself,
where the morphisms may only be pseudo. However, under
the present assumptions there is a left $2$-adjoint to the
inclusion of $\algebra{T}_\strict$ into $\algebra{T}$, sending an algebra $A$ to an
algebra $A'$ such that the category of strict morphisms
$A'\to B$ is isomorphic to the category of pseudo morphisms
$A\to B$. (See \cite[Theorem 3.13]{blackwell1989}.) Thus, when considering the category of diagrams
$d\colon J\to X$ in an object $X$ of $\algebra{T}$,
we can replace $d$ with its strict adjunct $\overline{d}\colon J'\to X$
and proceed using the results for $\algebra{T}_\strict$. As a practical matter, we will often have
specified $d$ to be strict in any case. Therefore the
results in the strict case suffice to classify the objects
of the diagram category even in the pseudo case.

% Bibliography

\printbibliography[heading=bibintoc]

\end{document}